

  \newcount\fontset
  \fontset=1

  \def\dualfont#1#2#3{\font#1=\ifnum\fontset=1 #2\else#3\fi}
  \dualfont\eightrm {cmr8} {cmr7}
  \dualfont\eightsl {cmsl8} {cmr7}
  \dualfont\eightit {cmti8} {cmti10}
  \dualfont\eightmi {cmmi8} {cmmi7}
  \dualfont\tensc {cmcsc10} {cmcsc10}
  \dualfont\eighttt {cmtt8} {cmtt10}
  \dualfont\titlefont {cmbx12} {cmbx10}
  \dualfont\eightsymbol {cmsy8} {cmsy10}


  \magnification=\magstep1
  \nopagenumbers
  \voffset=2\baselineskip
  \advance\vsize by -\voffset
  \headline{\ifnum\pageno=1 \hfil \else \tensc\hfil
    amenability for fell bundles
  \hfil\folio \fi}


  \def\vg#1{\ifx#1\null\null\else
    \ifx#1\ { }\else
    \ifx#1,,\else
    \ifx#1..\else
    \ifx#1;;\else
    \ifx#1::\else
    \ifx#1''\else
    \ifx#1--\else
    \ifx#1))\else
    { }#1\fi \fi \fi \fi \fi \fi \fi \fi \fi}
  \newcount\secno \secno=0
  \newcount\stno
  \def\goodbreak{\vskip0pt plus.1\vsize\penalty-250
    \vskip0pt plus-.1\vsize\bigskip}
  \outer\def\section#1{\stno=0
    \global\advance\secno by 1
    \goodbreak\vskip\parskip
    \message{\number\secno.#1}
    \noindent{\bf\number\secno.\enspace#1.\enspace}}
  \def\state#1 #2\par{\advance\stno by 1\medbreak\noindent
    {\bf\number\secno.\number\stno.\enspace #1.\enspace}{\sl #2}\medbreak}
  \def\nstate#1 #2#3\par{\state{#1} {#3}\par
    \edef#2{\number\secno.\number\stno}}
  \def\proof{\medbreak\noindent{\it Proof.\enspace}}
  \def\proofend{\ifmmode\eqno\square\else\hfill\square\medbreak\fi}
  \newcount\zitemno \zitemno=0
  \def\zitem{\global\advance\zitemno by 1 \smallskip
    \item{\ifcase\zitemno\or i\or ii\or iii\or iv\or v\or vi\or vii\or
    viii\or ix\or x\or xi\or xii\fi)}}
  \def\$#1{#1 $$$$ #1}


  \newcount\bibno \bibno=0
  \def\newbib#1{\advance\bibno by 1 \edef#1{\number\bibno}}
  \def\cite#1{{\rm[\bf #1\rm]}}
  \def\scite#1#2{{\rm[\bf #1\rm, #2]}}
  \def\loccit#1{(#1)}
  \def\se#1#2#3#4{\def\a{#1}\def\b{#2}\ifx\a\b#3\else#4\fi}
  \def\setem#1{\se{#1}{}{}{, #1}}
  \def\index#1{\smallskip \item{[#1]}}

  \def\ATarticle#1{\zarticle #1 xyzzy }
  \def\zarticle#1, author = #2,
   title = #3,
   journal = #4,
   year = #5,
   volume = #6,
   pages = #7,
   NULL#8 xyzzy {\index{#1} #2, ``#3'', {\sl #4\/} {\bf #6} (#5), #7.}

  \def\ATtechreport#1{\ztechreport #1 xyzzy }
  \def\ztechreport#1,
    author = #2,
    title = #3,
    institution = #4,
    year = #5,
    note = #6,
    toappear = #7,
    NULL#8
    xyzzy
    {\index{#1} #2, ``#3''\setem{#6}\setem{#4}\setem{#5}\se{#7}{}{.}{, to
appear in #7}}

  \def\ATbook#1{\zbook #1 xyzzy }
  \def\zbook#1,
    author = #2,
    title = #3,
    publisher = #4,
    year = #5,
    volume = #6,
    series = #7,
    NULL#8
    xyzzy
    {\index{#1} #2, ``#3''\setem{#7}\se{#6}{}{}{ vol. #6}, #4, #5.}


  \def\stress#1{{\it #1}\/}
  
  \def\crossproduct{\hbox to 1.8ex{$\times \kern-.45ex\vrule height1.1ex
depth0pt width0.45truept$\hfill}}
  \def\:{\colon}
  \def\*{\otimes}
  \def\+{\oplus}
  \def\x{\times}
  \def\({\bigl(}
  \def\){\bigl)}
  \def\|{\Vert}
  \def\<{\langle}
  \def\>{\rangle}
  \def\arw{\longrightarrow}
  \def\cstar{$C^*$}
  \def\square{\hbox{$\sqcap\!\!\!\!\sqcup$}}
  \def\for#1{,\quad #1}


  \def\arw{\rightarrow}
  \def\inv{^{-1}}
  \def\B{{\cal B}}
  \def\Lin{{\cal L}}
  \def\F{{\bf F}}
  \def\G{\Gamma}
  \def\L{\Lambda}
  \def\CB{C^*(\B)}
  \def\CrB{C^*_r(\B)}
  \def\CK{Cuntz--Krieger\vg}
  \def\OA{{\cal O}_A}
  \def\Fn{{\bf F}_n}
  \def\N{{\bf N}}
  \def\.{\cdot}
  \def\a{\alpha}
  \def\b{\beta}
  \def\g{\gamma}
  \def\s{\sigma}
  \def\e{{\rm e}}
  \def\claim#1{\medskip \noindent {\tensc CLAIM #1:}}


  \newbib\AEE
  \newbib\Claire
  \newbib\Blackadar
  \newbib\Bratteli
  \newbib\CKbib
  \newbib\SoftIII
  \newbib\SoftI
  \newbib\SoftII
  \newbib\Circle
  \newbib\AF
  \newbib\TPA
  \newbib\Unconditional
  \newbib\Inverse
  \newbib\FD
  \newbib\KK
  \newbib\McClanahan
  \newbib\Nica
  \newbib\Pedersen
  \newbib\Quigg
  \newbib\QR
  \newbib\Rieffel
  \newbib\SV
  \newbib\Wassermann
  \newbib\Woronowicz


  \null
  \vskip-2\bigskipamount

  \begingroup
  \def\c{\centerline}
  \c{{\titlefont AMENABILITY FOR FELL BUNDLES}\footnote{*}{\eightrm
  This paper is also available from {\eighttt
http://www.ime.usp.br/{\~{}}exel/}}}

  \bigskip
  \baselineskip=10pt
  \eightit
  \c{\tensc Ruy Exel\footnote{**}{\eightrm Partially supported by CNPq,
Brazil.}}
  \c{Departamento de Matem\'atica}

  \c{Universidade de S\~ao Paulo}
  \c{Rua do Mat\~ao, 1010}
  \c{05508-900 S\~ao Paulo -- Brazil}
  \c{exel@ime.usp.br}
  \endgroup


  \bigskip\bigskip
  \midinsert\narrower\narrower

  \noindent {\bf Abstract}.
  Given a Fell bundle $\B$, over a discrete group
  $\G$, we construct its reduced cross sectional algebra $\CrB$, in analogy
with the reduced crossed products defined for \cstar-dynamical systems.  When
the reduced and full cross sectional algebras of $\B$ are isomorphic, we say
that the bundle is amenable.  We then formulate an approximation property
which we prove to be a sufficient condition for amenability.

 A theory of $\G$-graded \cstar-algebras possessing a conditional expectation
is developed, with an eye on the Fell bundle that one naturally associates to
the grading.  We show, for instance, that all such algebras are isomorphic to
$\CrB$, when the bundle is amenable.

  We also study induced ideals in graded \cstar-algebras and obtain a
generalization of results of Str\v{a}til\v{a} and Voiculescu on AF-algebras,
and of Nica on quasi-lattice ordered groups.  A brief comment is made on the
relevance, to our theory, of a certain open problem in the theory of exact
\cstar-algebras.

  An application is given to the case of an $\Fn$--grading of the \CK algebras
$\OA$, recently discovered by Quigg and Raeburn.
Specifically, we show that the \CK bundle satisfies the
approximation property, and hence is amenable, for all matrices $A$ with
entries in $\{0,1\}$, even if $A$ does not satisfy the well known property (I)
studied by Cuntz and Krieger in their paper.
  \endinsert

  \section{Introduction}
  A possible definition of Fell bundles (also called \cstar-algebraic bundles
\cite{\FD}), for the special case of discrete groups, states that these are
given by a collection
  $\B = \(B_t\)_{t\in G}$
  of closed subspaces of a \cstar-algebra $B$, indexed by a discrete group
$\G$, satisfying
  $B_t^* = B_t$ and
  $B_t B_s \subseteq B_{ts}$ for all $t$ and $s$ in $\G$.
  If, in addition, the $B_t$'s are linearly independent and their direct sum
is dense in $B$, then $B$ is said to be a graded \cstar-algebra.

  Fell bundles and graded algebras occur in an increasing number of situations
in the theory of \cstar-algebras, often without the perception that they are
there.  The example in which the Fell bundle structure is the most
conspicuous, is that of the well known crossed product construction associated
to a \cstar-dynamical system \cite{\Pedersen}, recently extended to the case
of twisted partial dynamical systems \cite{\TPA}.  See also \cite{\Circle,
\McClanahan}.

In the most interesting cases, these bundles have a commutative unit fiber
algebra, that is, $B_e$ (where $e$ denotes the unit group element).  In that
case, by the main result of \cite{\TPA}, one can say that, up to
stabilization, the complexity of the bundle resides in three distinct
compartments, namely, the topological structure of the spectrum of $B_e$,
certain homeomorphisms between open sets in that spectrum, and a certain two
cocycle. See \cite{\TPA} for more details.

Among the examples in which the Fell bundle structure is not so striking to
the eye, lie some of the most intensely studied \cstar-algebras of the past
couple of decades.  These include all the
  AF-algebras \cite{\Bratteli,\SV,\AF},
  the \CK algebras $\OA$ \cite{\CKbib,\QR},
  algebras generated by Wiener-Hopf operators \cite{\Nica},
  non-commutative Heisenberg manifolds \cite{\Rieffel,\AEE},
  the quantum ${\rm SU}_2$ groups \cite{\Woronowicz} (here one may use one of
several easily available circle actions to obtain many interesting gradings,
as done in \cite{\Circle}),
  the soft tori $A_\varepsilon$ \cite{\SoftI, \SoftII, \SoftIII},
  and many others.

  However intriguing this may be, most of the examples cited possess a
commutative unit fiber algebra, and hence the comment above applies.
  But, there is a catch.  The bundle structure alone is not enough to
characterize the algebra.  The point is that non-isomorphic algebras may
possess identical associated Fell bundles, as in the case of the reduced and
full group \cstar-algebras of non-amenable discrete groups (see
\scite{\FD}{VIII.16.12}).

  It is the main purpose of this work to study this point in detail.  The crux
of the matter is thus to determine conditions on a Fell bundle $\B$, such that
all graded \cstar-algebras, whose associated Fell bundle coincides with $\B$,
are isomorphic to each other.  After we show that all such algebras lie in
between the full and reduced cross sectional algebras of the bundle, that is
$\CB$ and $\CrB$, respectively, this is equivalent to saying that the left
regular representation of $\CB$ is faithful.

Inspired by the work of Andu Nica \cite{\Nica}, we call such bundles
\stress{amenable}.
  Our main contribution is to formulate an approximation property for Fell
bundles, which we prove to be a sufficient condition for amenability.  This
condition is strongly influenced by the work of Claire Anantharaman-Delaroche
and, to a certain extent, could be thought of as an attempt at a
generalization of a similar condition studied in \cite{\Claire}.  We do not
claim to have taken the analogy to its limits, as the role of the center of
the unit fiber algebra, played in \cite{\Claire}, is yet to be understood in
the Fell bundle situation.

  Our major application is to the case of the recently discovered bundle
structure, over the free group $\Fn$, of the \CK algebras $\OA$, obtained by
Quigg and Raeburn in \cite{\QR}, in terms of a co-action of $\Fn$.  To study
this example, we show that the \CK relations, that is, the relations that
define the algebras $\OA$, give rise to a partial representation
\cite{\Inverse} of the free group.  By a partial representation of a group
$\G$ on a Hilbert space $H$, we mean a unital map $\s \: \G \arw \Lin(H)$,
  such that, for all $t,r\in \G$, one has
    $\s(t\inv) = \s(t)^*$, and
    $\s(t) \s(r) \s(r\inv) = \s(tr) \s(r\inv)$. See \cite{\Inverse} for more
details.

  In turn, given any partial representation of a discrete group, we construct
an associated Fell bundle. This construction ascribes the Fell bundle related
to $\OA$, mentioned above, to the partial representation arising from any
universal representation of the \CK relations.

  Irrespective of the main hypothesis imposed on the matrix $A$ in
\cite{\CKbib}, namely, condition (I), we show that the \CK bundle satisfies
the approximation property and hence is amenable.  In fact, we prove that this
holds for the Fell bundle associated to any partial representation $\s$ of the
free group, which satisfies
  $$
  \sum_{i=1}^n \s(g_i) \s(g_i)^* = 1,
  $$
  where $g_1,\ldots,g_n$ are the generators of the free group.

  On another front, we study induced ideals in graded algebras, where we are
able to mimic Nica's work in \cite{\Nica} and obtain results very similar to
his.  In fact, the work of Nica has been an inspiration for us all along, as
he also treats questions related to the approximation property.  Other than
merely formal generalizations, we seem to have gotten a little further,
because our proof of the amenability of the \CK bundle is given independently
of the work of Cuntz and Krieger, themselves, namely, the uniqueness of the
\cstar-algebra generated by non-trivial representations of the \CK relations,
when property (I) is present \cite{\CKbib}.

Still under the heading of induced ideals, we point out a curious relationship
with an open problem in the theory of exact \cstar-algebras, stated in
\scite{\Wassermann}{2.5.3}. The question of whether $C^*_r(\G)$ is an exact
\cstar-algebra for any countable discrete group $\G$, seems to be related to
what we do here.

  The author is indebted to several people who contributed in many ways for
the evolution of this work.  Among these he would like to express his thanks
to Claire Anantharaman-Delaroche for a brief, but fruitful discussion during
his short visit to Orleans, to Marcelo Laca for pointing out important
references and also for several interesting conversations, and to Cristina
Cerri for having carefully gone over the paper \cite{\Claire} at the operator
algebra seminar in S\~ao Paulo and for many discussions as well.

  \section {Reduced Cross Sectional Algebras}
  Throughout this section, $\G$ will denote a discrete group and
  $\B = \{B_t\}_{t\in G}$
  will be a fixed
  \cstar-algebraic bundle
  over $\G$, as defined in \cite{\FD}.  In recent years it has been customary
to refer to \cstar-algebraic bundles as \stress{Fell bundles}, a terminology
will think is quite appropriate.

Our specialization to the case of discrete groups has the primary purpose of
avoiding technical details which occur in the theory of Fell bundles over
continuous groups.  With some more work, we believe that much, if not all,
that we do here can be extended to the general case.  Another reason for us
being satisfied with the discrete group case is that our main application is
for $\G = \F_n$, the free group on $n$ generators.

By a \stress{section of} $\B$ we mean any function
  $$
  f \: \G \arw \bigcup_{t\in \G} B_t,
  $$
  such that
  $f(t) \in \B_t$ for all $t \in \G$.

Let
  $C_c(\B)$
  denote the set of all finitely supported sections of $\B$.  We shall regard
$C_c(\B)$ as a right module over $\B_e$, which, when equipped with the
$B_e$--valued inner product
  $$
  \<\xi,\eta\> = \sum_{t\in \G} \xi(t)^* \eta(t) \for \xi,\eta\in C_c(\B),
  $$
  becomes a pre-Hilbert module in the sense of \scite{\KK}{1.1.1}.  The
completion of $C_c(\B)$ can be shown to consist of all cross sections $\xi$ of
$\B$ such that the series
  $$
  \sum_{t\in \G} \xi(t)^* \xi(t)
  $$
  converges unconditionally.  Incidentally, we say that a series
  $\sum_{i\in I} x_i$
  in a Banach space $E$, indexed by any index set $I$, is unconditionally
summable to $x\in E$ if, for any
  $\varepsilon > 0$,
  there exists a finite subset
  $I_0 \subseteq I$
  such that, for all finite subsets
  $J \subseteq I$, with
  $I_0 \subseteq J$, one has
  $$
  \| x - \sum_{i\in J} x_i \| < \varepsilon.
  $$

As we said before, we believe in the possibility of generalizing the results
in this work to Fell bundles over continuous group.  If that is to be
attempted, then the unconditional summability, which is such a pervasive
ingredient here, is likely to be replaced by the concept of unconditional
integrability of \cite{\Unconditional}.

We denote by
  $l_2(\B)$
  the completion of $C_c(\B)$, so that $l_2(\B)$ becomes a right Hilbert
$\B_e$--module.  Given $\xi$ and $\eta$ in $l_2(\B)$ one can show that
  $$
  \<\xi,\eta\> = \sum_{t\in \G} \xi(t)^* \eta(t),
  $$
  where the series is unconditionally summable, as described above.

  As in \scite{\KK}{1.1.7}, we will denote by
  $\Lin_{B_e}(l_2(\B))$,
  or simply by $\Lin(l_2(\B))$, the \cstar-algebra of all adjointable
operators on $l_2(\B)$.

For each $t\in \G$ and each $b_t \in B_t$ (the subscript in $b_t$ is not
absolutely necessary but it will be used to remind us that $b_t$ belongs to
the fiber $B_t$), define
  $$
  \L(b_t) \xi |_s = b_t \xi(t\inv s)
  \for \xi\in l_2(\B)), \quad s\in \G.
  $$

  It is not difficult to show that
  $\L(b_t) \xi$ does belong to $l_2(\B)$, and that
  $\| \L(b_t) \xi \| \leq \| b_t \| \|\xi \|$.

  \state Proposition $\L(b_t)$ is an adjointable operator on $l_2(\B)$ and
$\L(b_t)^* = \L(b_t^*)$.

  \proof Given $\xi, \eta \in l_2(\B)$ we have
  $$
  \< \xi, \L(b_t) \eta \> =
  \sum_{s\in \G} \xi(s)^* b_t \eta(t\inv s) =
  \sum_{s\in \G} \xi(ts)^* b_t \eta(s) =
  \sum_{s\in \G} (b_t^* \xi(ts))^* \eta(s) = \ldots
  $$
  Observe that $b_t^* \in B_{t\inv}$, hence
  $ \L(b_t^*) \xi |_s = b_t^* \xi(t s). $
  So, the above equals
  $$
  \ldots =
  \sum_{s\in \G} (\L(b_t^*)\xi |_s)^* \eta(s) =
  \< \L(b_t^*) \xi, \eta\>.
  \proofend
  $$

  \nstate Proposition \LeftReg The map
  $$
  \L \: \bigcup_{t\in \G} B_t \arw \Lin(l_2(\B))
  $$
  is a representation of $\B$ in the sense that for all $t,s \in \G$ (compare
\scite{\FD}{VIII.9.1})
  \zitemno = 0
  \zitem $\L$ is a continuous linear map from $B_t$ to $\Lin(l_2(\B))$,
  \zitem for $b_t \in B_t$ and $c_s \in B_s$ one has $\L(b_t c_s) = \L(b_t)
\L(c_s)$,
  \zitem $\L(b_t^*) = \L(b_t)^*$.

  \proof Since (i) is easy and (iii) is already demonstrated, let us prove
(ii). Take $\xi \in l_2(\B)$, then, since $b_t c_s \in B_{ts}$,
  $$
  \L(b_t c_s) \xi |_r = b_t c_s \xi (s\inv t\inv r) =
  b_t\( \L(c_s) \xi |_{t\inv r} \) =
  \L(b_t) \L(c_s) \xi |_r.
  \proofend
  $$

  From now on we shall refer to $\L$ as the \stress{left regular
representation} of $\B$.

  \nstate Definition \Reduced The \stress{reduced cross sectional algebra} of
$\B$, denoted $\CrB$, is the sub \cstar-algebra of $\Lin(l_2(\B))$ generated
by the range of the left regular representation of $\B$.

We recall that the (full) cross sectional algebra of $\B$, denoted $\CB$
\scite{\FD}{VIII.17.2}, is the enveloping \cstar-algebra of $l_1(\B)$.  By the
universal property (see \scite{\FD}{VIII.16.12}) of $\CB$ one sees that there
is a canonical epimorphism of \cstar-algebras, which we call the
  \stress{left regular representation}
  of $\CB$ and, by abuse of language, still denote by $\L$,
  $$
  \L \: \CB \arw \CrB.
  $$

  We would now like to develop the rudiments of a Fourier analysis for $\CrB$,
and, in particular, to define Fourier coefficients for elements of $\CrB$. For
each $t\in \G$ let
  $$
  j_t \: B_t \arw l_2(\B)
  $$
  be the ``inclusion'' map, given by
  $$
  j_t(b_t)|_s = \left\{ \matrix{
    b_t &\hbox{if $s=t$} \cr
    0 &\hbox{if $s\neq t$} \cr
  }\right.
  $$
  for $b_t\in B_t$ and $s\in \G$

  We shall regard each $B_t$ as a right Hilbert $B_e$--module under the
obvious right module structure and the $B_e$--valued inner product given by
  $$
  \<b_t, c_t\> = b_t^* c_t \for b_t, c_t \in B_t.
  $$
  This said, one can easily prove that each $j_t$ is an adjointable map and
that for $\xi \in l_2(\B)$, one has
  $j_t^*(\xi) = \xi(t).$
  It follows that $j_t^* j_t$ is the identity map on $B_t$, and hence that
$j_t$ is an isometry.  This allows us to identify $B_t$ and its image
  $
  \bar B_t = j_t(B_t)
  $
  within $l_2(\B)$.

  Incidentally, this shows the very subtle fact that any adjointable map from
any Hilbert $B_e$--module into $l_2(\B)$, whose image lies in $\bar B_t$,
remains an adjointable map if its codomain is reduced to $\bar B_t$.

  \nstate Proposition \TwoFour Let $t,s \in \G$, $b_t \in B_t$ and $c_s \in
B_s$.  Then
  $\L(b_t) j_s(c_s) = j_{t s} (b_t c_s).$
  Therefore
  $\L(b_t) \bar B_s \subseteq \bar B_{t s}$.
  In addition, the map
  $$
  b_t \in B_t \mapsto \L(b_t) |_{\bar B_e} \in \Lin(\bar B_e, \bar B_t)
  $$
  is isometric.

  \proof We have
  $$\L(b_t) j_s(c_s)|_r =
  b_t\( j_s(c_s) |_{t\inv r}\) =
  \delta_{t\inv r, s} \ b_t c_s =
  j_{t s} (b_t c_s)|_r,
  $$
  where $\delta_{t\inv r, s}$ is the Kronecker symbol.  This proves the first
statement.  With respect to the isometric property of the map above, all we
are saying is that
  $$
  \sup\{ \| b_t a \| \: a \in B_e, \|a\| \leq 1 \} = \|b_t\|.
  $$
  This follows from the fact that any approximate unit for $B_e$ acts as an
approximate unit for $B_t$ also \scite{\FD}{VIII.16.3}.
  \proofend

  \nstate Corollary \Identify For each $t\in \G$ and each $b_t \in B_t$, one
has that
  $ \| \L(b_t) \| = \|b_t\| $
  and hence each $B_t$ may be identified with its image in $\CrB$.

  \proof Follows immediately from \loccit{\TwoFour}.
  \proofend

  The following is the crucial step in defining Fourier coefficients.

  \nstate Proposition \Fucoef For each $x$ in $\CrB$ and each $t$ in $\G$
there exists a unique $b_t$ in $B_t$ such that
  $j_t^* x j_e (a) = b_t a$ for every $a \in B_e$.  In addition we have
$\|b_t\| \leq \|x\|$.

  \proof Uniqueness follows from the fact that if $b_t a = 0$ for all $a \in
B_e$, then $b_t=0$.  As for the existence part, the easiest case, by far, is
when $B_e$ has a unit.  In this case $b_t$ is just
  $j_t^* x j_e(1)$.  If no unit is available suppose, to start with, that
  $x = \sum_{s\in \G} \L(b_s)$,
  where $b_s \in B_s$ and $b_s=0$ except for finitely many group elements $s$.
We then have
  $$
  j_t^* x j_e (a) = \sum_{s\in \G} j_t^* \L(b_s) j_e(a) = b_t a.
  $$
  This says that $j_t^* x j_e$, viewed as an element of $\Lin(B_e,B_t)$,
coincides with $\L(b_t) |_{B_e}$, and hence lies in the isometric copy of
$B_t$ within $\Lin(B_e,B_t)$ provided by \loccit{\TwoFour}.  Now, since the
set of all $x's$ considered is dense in $\CrB$, we have obtained the existence
part for all $x\in \CrB$.

  Next observe that
  $$
  \| b_t \| =
  \| \L(b_t) |_{B_e} \| =
  \| j_t^* x j_e\| \leq
  \|x \|.
  \proofend $$

  \state Definition For $x$ in $\CrB$ and $t\in \G$ the $t^{{\eightrm th}}$
  \stress{Fourier coefficient}
  of $x$ is the unique element $\hat x(t) \in B_t$ such that
  $j_t^* x j_e (a) = \hat x(t) a$
  for all $a$ in $B_e$.
  The
  \stress{Fourier transform}
  of $x$ is the cross section of $\B$ defined by $t \mapsto \hat x(t)$.

  Given the left regular representation of $\CB$
  $$
  \L \: \CB \arw \CrB,
  $$
  we can easily define the Fourier coefficients for elements of $\CB$ as well,
that is, if $y\in \CB$ we put $\hat y(t) = \widehat{\L(y)}(t)$.

  From the proof of \loccit{\Fucoef} we see that, if $x$ is the finite sum
  $x = \sum_{t\in \G} \L(b_t)$
  with $b_t \in B_t$, then $\hat x(t) = b_t$.  Also it may be worth insisting
that \loccit{\Fucoef} yields $\|\hat x(t) \| \leq \|x\|.$

  \nstate Proposition \DiagonalMatrix For $x$ in $\CrB$, $t,s \in \G$ and $b_t
\in B_t$ one has
  $$
  j_s^* x j_t(b_t) = \hat x(s t\inv) b_t
  $$

  \proof By continuity it is enough to consider finite sums
  $x = \sum_{r\in \G} \L(c_r)$
  as above.  For such an $x$ we have
  $$ j_s^* x j_t (b_t) =
  \sum_{r\in \G} j_s^* \L(c_r) j_t(b_t) =
  \sum_{r\in \G} j_s^* j_{rt}(c_r b_t) =
  c_{s t\inv} b_t =
  \hat x(s t\inv) b_t.
  $$
  We have used that $j_s^* j_r = 0$ when $s\neq r$, a fact that is easy to
see.
  \proofend

  The Fourier coefficient $\hat x(e)$ has special properties worth mentioning.

  \nstate Proposition \CondExp The map
  $E \: \CrB \arw B_e$ given by $E(x) = \hat x(e)$ is a positive, contractive
conditional expectation.

  \proof We first note that we are tacitly identifying $B_e$ and its sibling
$\L(B_e)$, as permitted by \loccit{\Identify}.  We have already seen that $E$
is a contractive map.  If $x = \L(b_e)$ with $b_e \in B_e$, then we saw that
$E(x)=b_e$, so $E$ is idempotent.

  If
  $x = \sum_{t\in \G} \L(b_t)$ as before, then
  $$
  x^* x = \sum_{t,s \in \G} \L(b_t)^* \L(b_s) =
  \sum_{t,s \in \G} \L(b_t^* b_s) =
  \sum_{t,r \in \G} \L(b_t^* b_{tr}) =
  \sum_{r\in \G} \L \( \sum_{t \in \G} b_t^* b_{tr} \).
  $$
  Now, since
  $\sum_{t \in \G} b_t^* b_{tr}$
  is in $B_r$, we have that
  $E(x^* x) = \sum_{t\in \G} b_t^* b_t \geq 0$,
  so $E$ is positive.  Finally, if $x = \sum_{t\in \G} \L(b_t)$ is a finite
sum, and if $a\in B_e$, then
  $$
  E(ax) = E\(\sum_{t\in \G} \L(ab_t)\) = ab_e = aE(x)
  $$
  and similarly, $E(xa) = E(x)a$.
  \proofend

  \state Proposition For $x\in \CrB$ one has that the sum
  $$
 \sum_{t\in \G} \hat x(t)^* \hat x(t)
  $$
  is unconditionally convergent, and hence
  $\xi_x = \(\hat x(t)\)_{t\in \G}$
  represents an element of $l_2(\B)$.  Also, for any $a$ in $B_e$ we have
  $ x j_e(a) = \xi_x a. $

  \proof Suppose that $x = \sum_{t\in \G} \L(b_t)$ is a finite sum, with $b_t
\in B_t$.  Then, obviously
  $\xi_x = \sum_{t\in \G}j_t(b_t)$
  belongs to $l_2(\B)$ and we have
  $$
  x j_e(a) = \sum_{t\in \G} \L(b_t)j_e(a) =
  \sum_{t\in \G} j_t(b_ta) =
  \sum_{t\in \G}j_t(b_t) a =
  \xi_x a.
  $$
  Also
  $$
  \| \xi_x\|^2 =
  \|\sum_{t\in \G} b_t^* b_t \| =
  \| E(x^* x) \| \leq
  \| x^* x\| = \|x\|^2.
  $$
  Therefore the relation
  $x \mapsto \xi_x \in l_2(\B)$
  defines a bounded map, which so far is defined only for the $x's$ as above,
but which may be extended to the whole of $\CrB$ by continuity.

  Again by continuity we have $xj_e(a) = \xi_x a$, for any $x$ in $\CrB$ and
any $a$ in $B_e$.  Next observe that, for $t\in \G$,
  $$
  j_t^* x j_e (a) =
  j_t^*(\xi_x a) =
  j_t^*(\xi_x) a.
  $$
  which implies, by \loccit{\Fucoef}, that
  $j_t^*(\xi_x) = \hat x(t)$.  This says that
  $\xi_x = \(\hat x(t)\)_{t\in \G}$
  and the proof is therefore concluded.
  \proofend

  \nstate Corollary \FormulaForE For $x\in \CrB$ one has
  $$
  E(x^*x) = \sum_{t\in \G} \hat x(t)^* \hat x(t).
  $$

  \proof For $a,b\in B_e$ we have
  $$
  a^*E(x^*x) b =
  \<a,\widehat{x^* x}(e)b\> =
  \<a,j_e^* x^* x j_e(b)\> =
  \<x j_e (a), x j_e(b)\>
  \$=
  \<\xi_x a, \xi_x b \> =
  a^* \<\xi_x, \xi_x \> b =
  a^* \sum_{t\in \G} \hat x(t)^* \hat x(t) b,
  $$
  from which the conclusion follows.
  \proofend

  With this we arrive at an important result, which says, basically that $E$
is a faithful conditional expectation on $\CrB$.

  \nstate Proposition \Faithful For $x\in \CrB$ the following are equivalent
  \zitemno = 0
  \zitem $E(x^* x) =0$
  \zitem $\hat x(t) = 0$ for every $t\in \G$
  \zitem $x=0$.

  \proof The equivalence of (i) and (ii) is a consequence of
\loccit{\FormulaForE}.  That (iii) implies (ii) is obvious, so let us prove
that (ii) implies (iii).  Assume $\hat x(t) = 0$ for all $t\in \G$.  Then, by
\loccit{\DiagonalMatrix}, it follows that $j_s^* x j_t = 0$ for all $t$ and
$s$ in $\G$.

  Now, note that any $\xi \in l_2(\B)$ is the sum of the unconditionally
convergent series
  $\xi = \sum_{t\in \G} j_t j_t^*(\xi)$.  So
  $$
  x\xi =
  \sum_{s\in \G} \sum_{t\in \G} j_s j_s^* x j_t j_t^* \xi = 0,
  $$
  that is, $x=0$.
  \proofend

  \section {Graded \cstar-algebras}
  In this section we will study the relationship between graded
\cstar-algebras and Fell bundles.  Graded algebras occur in a great number of
different contexts in the theory of operator algebras, as in the theory of
co-actions of discrete groups \cite{\Quigg}.  See also \cite{\TPA},
\cite{\Nica}, \cite{\Inverse}.

  The following concept is taken from \scite{\FD}{VIII.16.11}.

  \nstate Definition \Graded Let $B$ be a \cstar algebra, $\G$ be a discrete
group and let
  $\(B_t\)_{t\in \G}$
  be a collection of closed linear subspaces of $B$.  We say that
$\(B_t\)_{t\in \G}$ is a \stress{grading} for $B$ if, for each $t,s$ in $\G$
one has
  \zitemno=0
  \zitem $B_t^* = B_t$
  \zitem $B_t B_s \subseteq B_{ts}$
  \zitem The subspaces $B_t$ are independent and $B$ is the closure of the
direct sum
  $\bigoplus_{t\in \G} B_t.$

  \medskip \noindent In that case we say that $B$ is a $\G$--graded
\cstar-algebra.  Each $B_t$ is called a \stress{grading subspace}.

  The primary example of graded \cstar-algebras is offered by the theory of
Fell bundles.

  \state Proposition Let $\B$ be a Fell bundle over the discrete group
$\G$. Then $\CrB$ is a graded \cstar-algebra via the fibers $B_t$ of $\B$.

  \proof The only slightly non trivial axiom to be proved regards the
independence of the $B_t$'s.  Assume, for that purpose, that
  $x = \sum_{t\in \G} \L(b_t)$ is a finite sum with $b_t \in B_t$ and that
$x=0$.  Then $b_t = \hat x(t) = 0.$
  \proofend

  A similar reasoning shows that the full cross sectional algebra is also
naturally graded.
  Conversely, suppose we are given a graded \cstar-algebra
  $B = \overline{\bigoplus}_{t\in \G} B_t$
  (by this notation we wish to say that all of the conditions of
\loccit{\Graded} are verified).
  One can then construct a Fell bundle over $\G$, by taking the fibers of the
bundle to be the grading subspaces.  The multiplication and adjoint
operations, required on a Fell bundle, are defined by restricting the
corresponding operations on $B$.

  In fact there is a great similarity between the formal definitions of Fell
bundles and that of a graded \cstar-algebras.  However, there are important
conceptual differences, better illustrated by the example provided by the full
and reduced group \cstar-algebras of a non-amenable discrete group.  These are
non-isomorphic graded \cstar-algebras whose associated Fell bundles are
indistinguishable from each other. See also \scite{\FD}{VIII.16.11}.

  Suppose we are given a $\G$--graded \cstar-algebra
  $B = \overline{\bigoplus}_{t\in \G} B_t$.
  Let us denote by $\B$ its associated Fell bundle.  Then, by the universal
property of $\CB$ \scite{\FD}{VIII.16.11}, there is a unique \cstar-algebra
epimorphism
  $$
  \Phi \: \CB \arw B,
  $$
  which is the identity on each $B_t$ (identified both with a subspace of
$\CB$ and of $B$ in the natural way).  This says that $\CB$ is, in a sense,
the biggest \cstar-algebra whose associated Fell bundle is $\B$.  Our next
result will show that the reduced cross sectional algebra is on the other
extreme of the range.  It is also a very economical way to show a
\cstar-algebra to be graded.

  \nstate Theorem \Economical Let $B$ be a \cstar-algebra and assume that for
each $t$ in a discrete group $\G$, there is associated a closed linear
subspace
  $B_t \subseteq B$
  such that, for all $t$ and $s$ in $\G$ one has
  \zitemno=0
  \zitem $B_t B_s \subset B_{ts}$,
  \zitem $B_t^* = B_{t\inv}$,
  \zitem the closed linear span of \/ $\bigcup_{t\in \G} B_t$ is dense in $B$.

  \medskip \noindent Assume, in addition, that there is a bounded linear map
  $$
  F \: B \arw B_e,
  $$
  such that $F$ is the identity map on $B_e$ and that $F$ vanishes on each
$B_t$, for $t\neq e$.
  Then
  \smallskip \item{(a)} The subspaces $B_t$ are independent and hence
$\(B_t\)_{t\in \G}$ is a grading for $B$.
  \smallskip \item{(b)} $F$ is a positive, contractive, conditional
expectation.
  \smallskip \item{(c)} If $\B$ denotes the associated Fell bundle, then there
exists a \cstar-algebra epimorphism
  $$
  \lambda \: B \arw \CrB,
  $$
  such that
  $\lambda(b_t) = \L(b_t)$
  for each $t$ in $\G$ and each $b_t$ in $B_t$.

  \proof If $x=\sum_{t\in \G} b_t$ is a finite sum with $b_t \in B_t$, then
  $$
  x^*x =
  \sum_{t,s\in \G} b_t^* b_s =
  \sum_{r\in \G} \( \sum_{t\in \G} b_t^* b_{tr} \),
  $$
  so
  $F(x^*x) = \sum_{t\in \G} b_t^* b_t.$

  Therefore, if $x=0$ then each $b_t = 0$, which shows the independence of the
$B_t$'s, and also that $F$ is positive.  Given $a$ in $B_e$, it is easy to see
that $E(ax) = aE(x)$ and $E(xa) = E(x)a$.  So, apart from contractivity, (b)
is proven.

  Define a pre right Hilbert $B_e$--module structure on $B$ via the
$B_e$--valued inner product
  $$
  \<x,y\> = E(x^*y) \for x,y\in B.
  $$
  For $b,x \in B$ we have, by the positivity of $F$ that
  $$
  \<bx,bx\> = F(x^*b^*bx) \leq
  \|b\|^2 F(x^*x) =
  \|b\|^2 \<x,x\>.
  $$
  So, the left multiplication operators
  $$
  L_b \: x \in B \mapsto bx \in B
  $$
  are bounded and hence extend to the Hilbert module completion $X$ of $B$
(after moding out vectors of norm zero).
  It is then easy to show that
  $$
  L \: b \in B \mapsto L_b \in \Lin(X)
  $$
  is a \cstar-algebra homomorphism.

  Let
  $x=\sum_{t\in \G} b_t$ and
  $y=\sum_{t\in \G} c_t$
  be finite sums with $b_t,c_t\in B_t$ and regard both $x$ and $y$ as elements
of $X$.  We then have
  $$
  \<x,y\> =
  F(\sum_{t,s\in \G} b_t^* c_s) =
  F(\sum_{r\in \G} \sum_{t\in \G} b_t^* c_{tr})
  \$=
  \sum_{t\in \G} b_t^* c_t =
  \< \sum_{t\in \G} j_t(b_t), \sum_{t\in \G} j_t(c_t) \>,
  $$
  where the last inner product is that of $l_2(\B)$.  This is the key
ingredient in showing that the formula
  $$
  U(\sum_{t\in \G} b_t) = \sum_{t\in \G} j_t(b_t)
  $$
  can be used to define an isometry of Hilbert $B_e$--bimodules
  $$
  U \: X \arw l_2(\B).
  $$
  Here it is important to remark that the continuity of $F$ ensures that the
set of finite sums $\sum_{t\in \G} b_t$ is not only dense in $B$ but also in
$X$.

  For $b_t$ in $B_t$ and $c_t$ in $B_s$ we have
  $$
  U L_{b_t}(c_s) =
  U(b_t c_s) =
  j_{ts}(b_t c_s) =
  \L(b_t) j_s(c_s) =
  \L(b_t) U (c_s).
  $$
  Since the finite sums $\sum_{s\in \G} c_s$ are dense in $X$, as observed
above, we conclude that
  $$
  U L_{b_t} U^* = \L(b_t).
  $$
  This implies, for all $b$ in $B$, that the operator $U L_b U^*$ belongs to
$\CrB$. This defines a map
  $$
  \lambda \: b \in B \mapsto U L_b U^* \in \CrB,
  $$
  which satisfies the requirements in (c).

  The only remaining task, that is, the proof of contractivity of $F$, now
follows easily because $F = E \lambda$, where $E$ is the conditional
expectation of \loccit{\CondExp}.
  \proofend

  The map $\lambda$, above, should be thought of as a generalized left regular
representation of $B$.

{}From now on we will be mostly interested in graded algebras possessing a
conditional expectation, and so we make the following:

  \state Definition A grading
  $\(B_t\)_{t\in \G}$
  on the \cstar-algebra $B$ is said to be a \stress{topological grading} if
there exists a conditional expectation of $B$ onto $B_e$, as in
\loccit{\Economical}.

Recalling our discussion about $\CB$ being the biggest graded algebra for a
given Fell bundle, we now see that $\CrB$ is the smallest such, at least among
topologically graded algebras.

  \nstate Corollary \Ft Let
  $B = \overline{\bigoplus}_{t\in \G} B_t$
  be a topologically graded \cstar-algebra.  Then, for every $t$ in $\G$,
there exists a contractive linear map
  $$
  F_t \: B \arw B_t,
  $$
  such that, for all finite sums $x=\sum_{t\in \G} b_t$, with $b_t \in B_t$
one has
  $F_t(x) = b_t$.

  \proof Simply define
  $F_t(b) = \widehat{\lambda(b)} (t)$,
  where $\lambda$ is as in \loccit{\Economical}.
  \proofend

  \nstate Proposition \KerLRR Let
  $B = \overline{\bigoplus}_{t\in \G} B_t$
  be a topologically graded \cstar-algebra, with conditional expectation $F$,
and let $\lambda$ be its left regular representation, as in
\loccit{\Economical}. Then
  $$
  \hbox{Ker}(\lambda) =
  \{b\in B \: F(b^* b) =0 \}.
  $$

  \proof Observe, initially, that $F = E \lambda$ and hence
  $$
  F(b^* b) = E(\lambda(b)^* \lambda(b)),
  $$
  from where we see that $F(b^* b) = 0$ if and only if $\lambda(b) = 0$,
because $E$ is faithful, as proved in \loccit{\Faithful}.
  \proofend

  This can be employed to give a useful characterization of $\CrB$ among
graded algebras:

  \nstate Proposition \Charact Let
  $B = \overline{\bigoplus}_{t\in \G} B_t$
  be a topologically graded \cstar-algebra with faithful conditional
expectation.  Then $B$ is naturally isomorphic to $\CrB$.

  \proof The left regular representation
  $\lambda \: B \arw \CrB$
  of \loccit{\Economical}
  will be an isomorphism by \loccit{\KerLRR}.
  \proofend

  As an example, let us briefly treat the case of semi-direct product bundles,
as defined in \scite{\FD}{VIII.4} (see also \cite{\TPA}).  For this let
  $(A,\G,\alpha)$
  be a discrete \cstar-dynamical system, i.e., $A$ is a \cstar-algebra and
$\alpha$ is an action of the discrete group $\G$ on $A$.  The
\stress{semi-direct product bundle} associated to $(A,\G,\alpha)$ is defined
to be the Cartesian product $A \x \G$, with the operations
  \zitemno=0
  \zitem $z (a,t) + (b,t) = (z a + b, t)$
  \zitem $(a,r) (b,s) = (a \alpha_r(b), rs)$
  \zitem $(a,t)^* = (\alpha_t^{-1}(a), t\inv)$
  \zitem $\|(a,t)\| = \|a\|$
  \medskip \noindent
  for any $a,b\in A$, any $r,s,t\in \G$ and any complex number $z$.

  \state Proposition The reduced cross sectional algebra of the semi-direct
product bundle above is naturally isomorphic to the reduced crossed product
algebra
  $A \crossproduct_{\alpha, r} \G$.

  \proof $A \crossproduct_{\alpha, r} \G$ is graded and the conditional
expectation is faithful \scite{\Pedersen}{7.11.3}. The result then follows
{}from \loccit{\Charact}.
  \proofend

  Let us now conduct a study of induced ideals in graded algebras, inspired,
among others, by the work of Str\v{a}til\v{a} and Voiculescu on AF-algebras
\cite{\SV} and the work of Nica on quasi-lattice ordered groups \cite{\Nica}.

  \nstate Theorem \InducedI Let $J$ be a closed, two-sided ideal in a
topologically graded \cstar-algebra
  $B = \overline{\bigoplus}_{t\in \G} B_t$.  Denote by $F_t$ the projections
onto each $B_t$, provided by \loccit{\Ft}.  Define
  \zitemno=0
  \zitem $J_1=\left\langle J \cap B_e \right\rangle$, that is, the ideal
generated by $J \cap B_e$,
  \zitem $J_2=\left\{ b\in B: F_t(b) \in J, \forall t\in\G \right \}$,
  \zitem $J_3=\left\{ b\in B: F_e(b^* b) \in J \right \}$.
  \medskip\noindent Then $J_1 \subseteq J_2 = J_3$.

  \proof Let
  $\lambda \: B \arw \CrB$
  be as in \loccit{\Economical}, and recall from \loccit{\FormulaForE} that
  $$
  E(x^*x) = \sum_{t\in \G} \hat x(t)^* \hat x(t),
  $$
  for all $x$ in $\CrB$.  Now, since we have
  $F_t(b) = \widehat{\lambda(b)}(t)$, we conclude that
  $$
  F_e(b^*b) = \sum_{t\in \G} F_t(b)^* F_t(b) \for b\in B.
  $$
  Therefore, recalling that ideals are hereditary sub-algebras, we see that
  $F_e(b^*b) \in J$ if and only if $F_t(b)^* F_t(b) \in J$ for all $t$ in
$\G$, which is the same as saying that $F_t(b) \in J$ for all $t$.  This
proves that $J_2=J_3$.
  Now, observe that the ideal generated by
  $J \cap B_e$
  contains, as a dense set, the collection of finite sums of the form
  $$
  x = \sum_i b_{r_i} x_i c_{s_i},
  $$
  where
  $r_i,s_i\in \G$,
  $b_{r_i} \in B_{r_i}$,
  $c_{s_i} \in B_{s_i}$ and
  $x_i$ is in $J \cap B_e$.
  For $x$ as above, note that
  $$
  F_t(x) = \sum_{r_i s_i = t} b_{r_i} x_i c_{s_i},
  $$
  so $F_t(x) \subseteq J$.  This proves that $x\in J_2$ and, by the density of
the set of $x$'s considered, that $J_1\subseteq J_2$.
  \proofend

  Later, after concluding our study of amenable bundles, we shall return to
the question of the equality between $J_1$ and $J_2$.

  \state Definition We say that a closed, two-sided ideal $J$, in a graded
\cstar-algebra
  $B = \overline{\bigoplus}_{t\in \G} B_t$
  is an \stress{induced ideal} if
  $J = \left\langle J \cap B_e \right\rangle$.

  Before closing this section let us consider quotients of graded algebras.

  \nstate Proposition \Ideals Let $J$ be an induced ideal in a topologically
graded \cstar-algebra
  $B = \overline{\bigoplus}_{t\in \G} B_t$.  Then the quotient \cstar-algebra
$B/J$ is topologically graded by the spaces $\pi(B_t)$, where
  $\pi \: B \arw B/J$
  is the quotient map.

  \proof By \loccit{\InducedI}, $J$ is invariant under $F_t$.  Therefore $F_t$
gives a well defined bounded map on $B/J$, namely
  $$
  \tilde F_t (x+J) = F_t(x) + J \for x\in B.
  $$

  We claim that $\pi(B_t)$ is a closed subspace of $B/J$, for each $t$.  In
fact, let $x\in B_t$ and $y\in J$.  Then
  $$
  \| x - y \| \geq
  \| F_t(x-y) \| =
  \| x - F_t(y) \|.
  $$
  In view of the fact that $F_t(y) \in J \cap B_t$, we conclude that
  $$
  \| x + J \| =
  \inf \{ \| x - y \| \: y\in J \} =
  \inf \{ \| x - y \| \: y\in J \cap B_t \}.
  $$
  The last term above gives the norm, in
  $B_t / (J \cap B_t )$,
  of the element
  $x + (J \cap B_t)$.
  In other words, the natural map
  $$
  B_t / (J \cap B_t ) \arw B/J
  $$
  is an isometry.  Therefore,
  $\pi(B_t)$, being isometric to
  $B_t / (J \cap B_t )$,
  is a Banach space and hence closed in $B/J$.

  It is now immediate to verify that the collection
  $\( \pi(B_t) \)_{t\in \G}$
  satisfies \loccit{\Economical.i--iii}, and that $\tilde F_e$ fills in the
rest of the hypothesis there to allow us to conclude that
  $\( \pi(B_t) \)_{t\in \G}$ is, in fact, a grading for $B/J$.  That this is a
topological grading follows from \loccit{\Economical.b}.
  This concludes the proof.
  \proofend

  \section {Amenability for Fell bundles} We return to the case of a general
Fell bundle $\B$ over a discrete group $\G$.  Initially, we would like to
adopt a convention designed to simplify our notation.  That is, we will
enforce, to the fullest extent, the identification of the fibers $B_t$ of $\B$
with its isomorphic image in any $\G$--graded \cstar-algebra in sight, whose
associated Fell bundle is identical to $\B$.  Of course, this applies to
$\CrB$ and $\CB$.  An excellent point of view is, in fact, to regard all such
graded \cstar-algebras as completions of the *-algebra
  $\bigoplus_{t\in \G} B_t$
  under various norms.

  The result in \loccit{\Economical}, together with the remark preceding it,
tells us that the topologically graded algebras whose associated Fell bundles
coincide with a given Fell bundle $\B$, are to be found among the quotients of
$\CB$ by ideals contained in the kernel of the left regular representation
  $\L \: \CB \arw \CrB$
  mentioned after \loccit{\Reduced}.

  It is therefore crucial to understand the kernel of $\L$ and, in particular,
to determine conditions to ensure that $\L$ is injective.  While we believe a
full comprehension of the problem is beyond the present methods, we intend to
give a sufficient condition for the injectivity of $\L$, which is general
enough to be applicable to an interesting Fell bundle appearing in the theory
of \CK algebras.  First some terminology.

  \state Definition The Fell bundle $\B$ is said to be \stress{amenable} if
the left regular representation
  $$
  \L \: \CB \arw \CrB.
  $$
  is injective.

  An immediate consequence of the results in the previous section is:

  \nstate Proposition \Democracy Let $\B$ be an amenable Fell bundle.  Then
all topologically graded \cstar-algebras, whose associated Fell bundles
coincide with $\B$, are isomorphic to each other.

  One of the main tools in the technique employed below is a study of
functions
  $a\: \G \arw B_e$
  such that
  $\sum_{t\in \G} a_t^* a_t$ converges unconditionally.  The reason is that it
allows us to construct ``wrong way maps''
  $\Psi \: \CrB \arw \CB$, as we will now see.

  \nstate Lemma \WrongWay Let
  $a\: \G \arw B_e$
  be such that
  $\sum_{t\in \G} a(t)^* a(t)$
  converges unconditionally.  Then there exists a bounded, completely positive
map
  $$
  \Psi \: \CrB \arw \CB
  $$
  such that
  $$
  \Psi(b_t) = \sum_{r\in \G} a(tr)^* b_t a(r) \for b_t\in B_t.
  $$
  In addition
  $\| \Psi \| \leq \| \sum_{t\in \G} a(t)^* a(t) \|.$

  \proof Assume, as we may, that $\CB$ is faithfully represented as an algebra
of operators on a Hilbert space $H$.  In fact, in order to lighten up our
notation, we will think of $\CB$ as a subalgebra of $\Lin(H)$.  This also
gives us a representation of $\B$ on $H$, as in \scite{\FD}{VIII.9.1}.  Define
a new representation of $\B$ on
  $H\* l_2(\G)$
  by letting
  $$
  \pi(b_t) = b_t \* \lambda_t \for t\in \G, \quad b_t\in B_t,
  $$
  where $\lambda \: \G \arw \Lin(l_2(\G))$ is the usual left regular
representation of $\G$.  Let $B$ be the \cstar-algebra of operators on
  $H\* l_2(\G)$
  generated by all the $\pi(b_t)$.

  Denoting the canonical basis of $l_2(\G)$ by
  $\(\delta_t\)_{t\in \G}$
  consider the isometries
  $$
  j_t \: \xi \in H \arw \xi\* \delta_t \in H \* l_2(\G).
  $$

  For each finite sum
  $ b = \sum_{t\in \G} \pi(b_t) $,
  with $b_t\in B_t$, one can easily verify that
  $j_t^* b j_s = b_{t s\inv}$.  In particular, this shows that $b=0$ causes
each $b_t$ to vanish.  It follows that
  $\( B_t \* \lambda_t \)_{t\in \G}$
  is a grading for $B$.  Moreover, we claim that the natural conditional
expectation
  $$
  F \: b \in B \arw (j_e^* b j_e) \* \lambda_e \in B_e \* \lambda_e
  $$
  is faithful.  In fact suppose that $F(b^*b)=0$, for a given $b$ in $B$.
Then $j_e^* b^* b j_e = 0$ and hence $b j_e = 0$.  Consider the blown up right
regular representation of $\G$, that is, the (anti) representation of $\G$ on
$H\* l_2(\G)$ given by
  $$
  \rho_t(\xi \* \delta_s) = \xi \* \delta_{st}
  \for s,t\in \G, \quad \xi \in H.
  $$
  It is easy to see that each $\rho_t$ belongs to the commutant of $B$ and
also that
  $$
  \rho_t j_s = j_{st} \for s,t \in \G.
  $$
  Now, since we have
  $b j_e = 0$, we will also have
  $b j_t = b \rho_t j_e = \rho_t b j_e = 0$.
  This implies that $b$ vanishes on each $H\*\delta_t$ and hence that $b=0$.
  With our claim proven we may apply \loccit{\Charact} to conclude that $B$ is
naturally isomorphic to $\CrB$.

  Consider the operator
  $$
  V \: H \arw H \* l_2(\G)
  $$
  given by the formula
  $V(\xi) = \sum_{t\in \G} a(t) \xi \* \delta_t,$
  for $\xi \in H.$
  This is well defined because
  $$
  \sum_{t\in \G} \| a(t) \xi \|^2 =
  \sum_{t\in \G} \< a(t)^* a(t) \xi, \xi \> \leq
  \| \sum_{t\in \G} a(t)^* a(t) \| \| \xi \|^2,
  $$
  which also informs us that $V$ is bounded with
  $\| V \| \leq \| \sum_{t\in \G} a(t)^* a(t) \|^{1/2}$.

  Define the completely positive map
  $$
  \Psi \: \Lin(H\* l_2(\G)) \arw \Lin(H)
  $$
  by $\Psi(x) = V^* x V$, for each $x$ in $\Lin(H\* l_2(\G))$.

  Observe that for every $b_t$ in $B_t$ and $\xi$ in $H$ one has
  $$
  \Psi(b_t \* \lambda_t) \xi =
  V^* (b_t \* \lambda_t) \( \sum_{r\in \G} a(r) \xi \* \delta_r \) =
  V^* \( \sum_{r\in \G} b_t a(r) \xi \* \delta_{tr} \) =
  \sum_{r\in \G} a(tr)^* b_t a(r) \xi,
  $$
  where we have used the fact that
  $$
  V^*(\eta \* \delta_t) =
  a(t)^* \eta \for t\in \G, \quad \eta\in H,
  $$
  which the reader will have no difficulty in proving.
  Summarizing, have that
  $$
  \Psi(b_t \* \lambda_t) = \sum_{r\in \G} a(tr)^* b_t a(r).
  $$

  Recall that $B$ is generated by the $b_t\*\lambda_t$. So, we see that $\Psi$
maps $B$ into $\CB$.  Therefore, if we identify $\CrB$ with $B$, we may
restrict $\Psi$ to $\CrB$ and get the desired map.

  Finally
  $$
  \| \Psi \| \leq
  \| V \|^2 \leq
  \| \sum_{t\in \G} a(t)^* a(t) \|.
  \proofend
  $$

  This puts us in position to describe an important concept.

  \nstate Definition \AP We say that the Fell bundle $\B$ has the
\stress{approximation property} if there exists a net
  $\( a_i \)_{i\in I}$
  of functions
  $a_i\: \G \arw B_e$,
  which is uniformly bounded in the sense that
  $$
  \sup_{i\in I} \| \sum_{t\in \G} a_i(t)^* a_i(t) \| < \infty,
  $$
  and such that for all $b_t$ in each $B_t$ one has
  $$
  \lim_{i \rightarrow \infty} \sum_{r\in \G} a_i(tr)^* b_t a_i(r)=
  b_t.
  $$

  An equivalent, although apparently stronger form of the approximation
property is given by our next result:

  \nstate Proposition \Stronger Suppose $\B$ has the approximation property.
Then there exists a net
  $\( a_i \)_{i\in I}$
  of \stress{finitely supported} functions
  $a_i\: \G \arw B_e$,
  satisfying the conditions of \loccit{\AP}.

  \proof The approximation property implies that there exists a constant $M>0$
such that for all finite sets
  $$
  X=\{b_{t_1},b_{t_2}, \ldots, b_{t_n} \} \subseteq \bigcup_{t\in\G}B_t
  $$
  and any $\varepsilon > 0$, there exists a function
  $a\: \G \arw B_e$
  such that
  $\| \sum_{t\in \G} a(t)^* a(t) \| < M$
  and
  $
  \| \sum_{r\in \G} a(t_k r)^* b_{t_k} a(r) \| < \varepsilon \for
k=1,\ldots,n.
  $
  By taking the restriction of $a$ to a sufficiently large set, we may assume
that $a$ is finitely supported.  Once this is done, denote such an $a$ by
$a_{X,\varepsilon}$ and think of it as a net, indexed by the set of pairs
$(X,\varepsilon)$, where $X$ is any finite subset of
  $\bigcup_{t\in\G}B_t$ and $\varepsilon$ is a positive real.
  It is often the case that the best nets are indexed by the most awful
looking sets.

  If these pairs are ordered by saying that
  $(X_1,\varepsilon_1) \leq (X_2,\varepsilon_2)$
  if and only if
  $X_1 \subseteq X_2$ and $\varepsilon_1 \geq \varepsilon_2$
  then the net
  $\( a_{X,\varepsilon}\)_{(X,\varepsilon)}$
  satisfies the required properties.
  \proofend

  The relationship between the approximation property and the amenability of
Fell bundles is one of our main results:

  \nstate Theorem \Main If a Fell bundle $\B$ has the approximation property,
then it is amenable.

  \proof Let
  $\( a_i \)_{i\in I}$
  be as above. For each $i\in I$ consider the map
  $\Psi_i \: \CrB \arw \CB$
  provided by \loccit{\WrongWay}, regarding each $a_i$.
  Define
  $\Phi_i \: \CB \arw \CB$
  to be the composition
  $\Phi_i = \Psi_i \Lambda$
  and observe that, by hypothesis
  $\lim_{i \rightarrow \infty} \Phi_i(b_t) = b_t$, for every $b_t$ in each
$B_t$.  Because the $b_t$'s span a dense subspace of $\CB$, and because the
$\Phi_i$'s are uniformly bounded, we conclude that
  $\lim_{i \rightarrow \infty} \Phi_i(x) = x$
  for every $x$ in $\CB$.

  Now, if $x\in \CB$ is in the kernel of the left regular representation, that
is, if $\L(x)= 0$, then
  $$
  x =
  \lim_{i \rightarrow \infty} \Phi_i(x) =
  \lim_{i \rightarrow \infty} \Psi_i(\L(x)) = 0,
  $$
  which proves that $\L$ is injective.
  \proofend

  A somewhat trivial example of this situation is that of bundles over
amenable groups.

  \state Theorem Let $\B$ be a Fell bundle over the discrete amenable group
$\G$. Then $\B$ satisfies the approximation property and hence is amenable.

  \proof According to \scite{\Pedersen}{7.3.8}, there exists a bounded net
  $\(f_i\)_{i\in I} \subseteq l_2(\G)$
  such that
  $$
  \lim_{i \rightarrow \infty} \sum_{r\in \G} \overline{ f_i(tr) } f_i(r) =
  1 \for t\in \G.
  $$
  Let
  $\(u_j\)_{j\in J}$
  be an approximate unit for $B_e$.  One then checks that the doubly indexed
net
  $\( a_{i,j} \)_{(i,j)\in I\x J}$
  defined by
  $a_{i,j}(t) = f_i(t) u_j$
  satisfies the definition of the approximation property.
  \proofend

  A very important consequence of the approximation property is the existence
of an analogue of Cesaro sums, which we would now like to discuss.

  \nstate Definition \SumProc A bounded linear map $\Phi \:B \arw B$, where
  $B = \overline{\bigoplus}_{t\in \G} B_t$
  is a graded \cstar-algebra, is said to be a \stress{summation process} if,
for all $x$ in $B$
  \zitemno = 0
  \zitem $\sum_{t\in \G} \Phi(\hat x(t))$ is unconditionally convergent,
  \zitem $\Phi(x) = \sum_{t\in \G} \Phi(\hat x(t))$,
  \zitem $\Phi(B_t) \subseteq B_t \for t\in \G$.

  \nstate Proposition \Cesaro Let
  $B = \overline{\bigoplus}_{t\in \G} B_t$
  be a topologically graded \cstar-algebra.  Assume that the Fell bundle $\B$,
associated to $B$, has the approximation property.  Then there exists a
bounded net $\(\Phi_i\)_{i\in I}$ of summation processes, converging to the
identity map of $B$ in the pointwise topology.

  \proof Since $\B$ is amenable, by \loccit{\Democracy}, $B$ is isomorphic to
$\CrB$.  We may then assume that, in fact, $B$ is the reduced cross sectional
algebra $\CrB$, itself.

  Let $\( a_i \)_{i\in I}$ be a finitely supported net, as in
\loccit{\Stronger}.  Denote by $\Psi_i$ the completely positive map
  $\Psi \: \CrB \arw \CB$
  given by \loccit{\WrongWay}, regarding each $a_i$, and let
  $\Phi_i = \Psi_i \Lambda$, following the first steps of the proof of
\loccit{\Main}. So, we have
  $$
  \lim_{i \rightarrow \infty} \Phi_i(x) = x \for x \in \CB.
  $$

  For each $i\in I$, let $S_i$ be the finite support of $a_i$, so
  $$
  \Phi_i(b_t) =
  \sum_{r\in \G} a_i(tr)^* b_t a_i(r) =
  \sum_{r \in S_i} a_i(tr)^* b_t a_i(r) \for b_t\in B_t.
  $$
  One can see, from this, that $\Phi_i(b_t) = 0$ unless $t\in S_i
S_i\inv$. This says that the series appearing in \loccit{\SumProc.i} is in
fact a finite series, hence convergent.

  Also, \loccit{\SumProc.iii} follows from the above formula for $\Phi_i$,
since each $a_i(r) \in B_e$.

  Finally, if $x$ is a finite sum
  $x=\sum_{t\in \G} b_t$, with
  $b_t \in B_t$, then $\hat x(t)=b_t$, that is
  $x=\sum_{t\in \G} \hat x(t)$.  This implies \loccit{\SumProc.ii} for such an
$x$ and hence also for all $x\in B$, by continuity of both sides of
\loccit{\SumProc.ii}, now that we now that the sum involved is finite.
  \proofend

Recalling our discussion of induced ideals in \loccit{\InducedI} we may now
add the following:

  \nstate Proposition \InducedII Let
  $B = \overline{\bigoplus}_{t\in \G} B_t$
  be a topologically graded \cstar-algebra.  Assume, in addition to the
hypothesis of \loccit{\InducedI}, that the Fell bundle $\B$, associated to
$B$, has the approximation property.  Then
  $J_1 = J_2 = J_3$.

  \proof Let $\(\Phi_i\)_{i\in I}$ be a net of summation processes, converging
pointwise to the identity map of $B$, as in \loccit{\Cesaro}.

  Next, let $x$ be in $J_2$, so that $\hat x(t) \in J \cap B_t$ for all $t\in
\G$ and, consequently $\hat x(t)^* \hat x(t)\in J \cap B_e$.  Recall that any
element $x$, in any \cstar-algebra, satisfies
  $x = \lim_{n\rightarrow \infty} x(x^*x)^{1/n}.$
  This implies that our $\hat x(t)$ belongs to the ideal generated by $J \cap
B_e$, that is, $J_1$.

  On the other hand, since
  $$
  \Phi_i(\hat x(t)) = \sum_{r\in \G} a_i(tr)^* \hat x(t) a(r),
  $$
  we have that
  $\Phi_i(\hat x(t))$ is in $J_1$, and that
  $\Phi_i(x) = \sum_{t\in \G} \Phi_i(\hat x(t))$
  is also in $J_1$.
  Finally, since
  $x = \lim_{i \rightarrow \infty} \Phi_i(x)$,
  we have that $x\in J_1$.
  \proofend

The result above can definitely fail in the absence of the approximation
property.  A counter-example is provided by the zero ideal in the full group
\cstar-algebra of a non-amenable group.  In this case $J_1 = \{0\}$ while
$J_3$ is the non-trivial kernel of the left regular representation, by
\loccit{\KerLRR}.

It is then natural to ask whether the result of \loccit{\InducedII} holds,
without the approximation property, but for the special case of the reduced
cross sectional algebra $\CrB$.  This turns out to be a very delicate
question, relating to the theory recently developed by E.~Kirchberg of exact
\cstar-algebras.  See \cite{\Wassermann} for details.  It still seems to be
unknown whether the reduced group \cstar-algebra of a countable discrete group
is exact \scite{\Wassermann}{2.5.3}.

Suppose, for a moment, that there exists a discrete non-amenable group $\G$,
for which $C^*_r(\G)$ is not exact and let, therefore
  $$
  0 \arw I \arw A \arw B \arw 0
  $$
  be an exact sequence of \cstar-algebras for which
  $$
  0 \arw I \* C^*_r(\G) \arw A \* C^*_r(\G) \arw B \* C^*_r(\G) \arw 0
  $$
  is not exact, where we are taking the spatial tensor product.

Consider the trivial group action of $\G$ on $A$ and recall that the reduced
cross sectional algebra of the semi-direct product bundle $(A,\G,id)$, i.e.,
the reduced crossed product
  $A \crossproduct_{id, r} \G$,
  coincides with the spatial tensor product $A \* C^*_r(\G)$.  Let $J$ be the
image, in
  $A \* C^*_r(\G)$,
  of
  $I \* C^*_r(\G)$.
  Referring to the sets $J_1$, $J_2$ and $J_3$ of \loccit{\InducedI}, one may
now prove that
  $J_1 = J = I \* C^*_r(\G)$,
  while $J_2$ is the kernel of the map
  $A \* C^*_r(\G) \arw B \* C^*_r(\G)$,
  which is strictly bigger than $J_1$, by hypothesis.

We would now like to explore an example of a Fell bundle over the non-amenable
free group which will be shown to satisfy the approximation property.  The
last two sections of the present work will be devoted to this.

  \section {The \CK bundle} The main goal of this section will be to describe
a Fell bundle over the free group $\Fn$ on $n$ generators, which arises from a
grading of the \CK algebras $\OA$.  This grading was discovered by Quigg and
Raeburn \cite{\QR} in the form of a co-action of the free group on $\OA$.

We will, however, prefer to give a description of this bundle in terms of
partial group representations because, on the one hand, this is a framework in
which a significant generalization is possible and, on the other, the
algebraic manipulations seem especially well suited to treat the problem at
hand.

  Let
  $A = \(a_{i,j}\)_{i,j}$
  be an $n\x n$ matrix in which
  $a_{i,j} \in \{0,1\}$,
  and suppose, furthermore, that no row or column of $A$ is identically
zero. In \cite{\CKbib} Cuntz and Krieger studied $n$--tuples
  $(S_1,S_2,\ldots,S_n)$
  of operators on a Hilbert space $H$ satisfying what we shall call the \CK
relations, namely
  \smallskip \itemitem{(CK0)}
  $S_i S_i^* S_i = S_i$, that is, each $S_i$ is a partial isometry,
  \smallskip \itemitem{(CK1)} $S_i^* S_j = 0$, for $i\neq j$,
  \smallskip \itemitem{(CK2)} $\sum_{i=1}^n S_i S_i^* =1$,
  \smallskip \itemitem{(CK3)} $S_i^* S_i = \sum_{j=1}^n a_{i,j} S_j S_j^*$.

  \medskip
  The \CK algebra $\OA$ can be defined as being the universal
\cite{\Blackadar} \cstar-algebra generated by $n$ symbols,
  $S_1,S_2,\ldots,S_n$,
  subject to the \CK relations.  One of the main results in \cite{\CKbib} is
that, when $A$ satisfies a certain condition, referred to as condition (I) in
\cite{\CKbib}, then any $n$--tuple of nonzero operators
  $(S_1,S_2,\ldots,S_n)$,
  satisfying \loccit{CK0--4}, generate a \cstar-algebra that is isomorphic of
$\OA$.

  In our treatment of $\OA$ we will make no special requirements on $A$ and so
we must insist on the definition of $\OA$ as a universal \cstar-algebra,
observing that, in the absence of condition (I), one may find \cstar-algebras
generated by partial isometries satisfying the above relations which are not
isomorphic to $\OA$.

The \CK relations are intimately related to the concept of partial
representations for $\Fn$, as we would like to show.  Recall from
\cite{\Inverse}, that a partial representation of a group $\G$ on a Hilbert
space $H$ is a map
  $S \: \G \arw \Lin(H)$
  satisfying, for all $t,r\in \G$,
  \smallskip \itemitem{(PR1)}
    $S(e) = I$, the identity operator on $H$,
  \smallskip \itemitem{(PR2)}
    $S(t\inv) = S(t)^*$, and
  \smallskip \itemitem{(PR3)}
    $S(t) S(r) S(r\inv) = S(tr) S(r\inv)$.

  \medskip
  Let us fix, for the time being, an $n$--tuple
  $(S_1,S_2,\ldots,S_n)$
  of operators on $\Lin(H)$ satisfying \loccit{CK0--4}.  We will construct,
based on these, a partial representation of $\Fn$.  For this, let
  $G = \{g_1,g_2,\ldots,g_n\}$
  be a set of (free) generators of $\Fn$.  It is well known that every element
$t$ of $\Fn$ has a unique decomposition (called its reduced decomposition, or
reduced form)
  $$
  t = x_1 x_2 \ldots x_k,
  $$
  where $x_i\in G \cup G\inv$ and $x_{i+1} \neq x_i\inv$ for all $i$.

  In this case, the integer $k$ will be referred to as the \stress{length} of
$t$ and we will write $|t|=k$.  It is elementary to show that $|rs| \leq |r| +
|s|$.

  If $t = x_1 x_2 \ldots x_k$ is as above, define
  $$
  S(t) = S(x_1) S(x_2) \ldots S(x_k),
  $$
  where we put
  $S(x) = S_j$ if $x=g_j$ or
  $S(x) = S_j^*$ if $x=g_j\inv$.

  One should not expect $S$ to be a multiplicative map, that is $S(tr)$ will
often differ from $S(t)S(r)$.  A significant case in which these coincide is
when the reduced forms of $t$ and $r$
  $$
  t = x_1 x_2 \ldots x_k
  $$
  $$
  r = y_1 y_2 \ldots y_k
  $$
  satisfy $x_k \neq y_1\inv$ because, then, the reduced form of $tr$ will be
  $t = x_1 x_2 \ldots x_k y_1 y_2 \ldots y_k$
  and we will have
  $S(tr) = S(t)S(r)$, by definition.

Note that the above condition on the pair $(t,r)$ is equivalent to saying that
$|tr|=|t|+|r|$.  In fact, this property will acquire some relevance below, so
we believe it justifies a special notation.

  \state Definition If $t,r\in \Fn$ are such that $|tr|=|t|+|r|$, then the
product $tr$ will be denoted by $t\.r$.

  Of course $t\.r$ is not supposed to indicate any new operation on $\Fn$ as
it means nothing but the usual product of $t$ and $r$.  However, it reminds us
that no cancelation is taking place when the reduced forms of $t$ and $r$ are
written side by side, in that order.

The first relevant use of our new notation is, thus, in the formula
  $$
  S(t\.r) = S(t) S(r)
  $$
  which holds whenever $|tr|=|t|+|r|$.

  Let $P$ be the subset of $\Fn$ consisting of elements
  $\alpha = x_1 x_2 \ldots x_k$
  where each $x_i$ is in $G$ (as opposed to $G\cup G\inv$). Such elements will
be called \stress{positive} and will usually be denoted by letters taken from
the beginning of the Greek alphabet.  For each natural number $k$ we will
denote by $P_k$ the set of positive elements of length $k$.  By convention,
$P_0 = \{e\}$, the singleton containing the identity element of $\Fn$.

  If $t = x_1 x_2 \ldots x_k$ is in reduced form and if, for some $i$, we have
$x_i \in G\inv$ and $x_{i+1}\in G$ then, by \loccit{CK1} we will have $S(x_i)
S(x_{i+1}) = 0$ and hence $S(t)=0$.  Therefore, the only case when $S(t)$ has
a chance of not vanishing is when $t=\a\b\inv$ for some $\a,\b\in P$ (although
$S(t)$ could vanish even in this case).

  \nstate Theorem \PR The map $S\: t\in \F_n \arw S(t) \in \Lin(H)$ is a
partial representation of \/ $\Fn$.

  \proof Properties \loccit{PR1--2} are trivial, so we concentrate on
\loccit{PR3}.  Our proof will consist of a series of claims, the first of
which is:
  \claim 1 For $\a\in P$ and $g_j\in G$ one has
  $S(\a g_j )^* S(\a g_j ) = \varepsilon S(g_j )^* S(g_j )$
  where $\varepsilon$ is either zero or one.

  \medskip
  This is obvious if $|\a|=0$ so, using induction, assume
  $\a=\b g_i $ with
  $\b\in P$ and
  $g_i \in G$. Then
  $$
  S(\a g_j )^* S(\a g_j ) =
  S(\b g_i \. g_j )^* S(\b g_i \. g_j ) =
  S(g_j)^* S(\b g_i)^* S(\b g_i) S(g_j)
  \$=
  \varepsilon S(g_j)^* S(g_i)^* S(g_i) S(g_j) =
  \varepsilon S(g_j)^* \sum_{l=1}^n a_{ik} S(g_k) S(g_k)^* S(g_j)
  \$=
  \varepsilon a_{ij} S(g_j)^* S(g_j) S(g_j)^* S(g_j) =
  \varepsilon a_{ij} S(g_j)^* S(g_j)
  $$
  proving our first claim.

  This shows, in particular, that for $\b\in P$ one has that $S(\b)^*S(\b)$ is
an idempotent.  Now, every operator $T$, such that $T^*T$ is idempotent, must
be a partial isometric operator and, thus $S(\b)$ is a partial isometry.

  \claim 2 For every $\a$ and $\b$ in $P$, if $|\a|=|\b|$ and $S(\a)^*S(\b)
\neq 0$ then $\a = \b$.

  \medskip Let $m=|\a|=|\b|$.  If $m=1$ then the claim is a consequence of
\loccit{CK1}. If $m \geq 1$ write
  $\a = \tilde \a g_i$ and
  $\b = \tilde \b g_j$ with
  $\tilde \a, \tilde \b \in P$ and
  $g_i,g_j \in G$.  Then
  $$
  0 \neq
  S(\tilde \a \. g_i)^* S(\tilde \b \. g_j) =
  S(g_i)^* S(\tilde \a)^* S(\tilde \b) S(g_j).
  $$
  So, in particular
  $ S(\tilde \a)^* S(\tilde \b) \neq 0$
  and, by induction, we have
  $\tilde \a = \tilde \b$. By claim (1) it follows that
  $S(\tilde \a)^* S(\tilde \a) = \varepsilon S(g_k )^* S(g_k)$
  for some $g_k\in G$.  Therefore
  $$
  S(\a)^* S(\b) =
  \varepsilon S(g_i)^* S(g_k)^* S(g_k) S(g_j)
  \$=
  \varepsilon S(g_i)^* \sum_{l=1}^n a_{kl} S(g_l) S(g_l)^* S(g_j) =
  \varepsilon a_{kj} S(g_i)^* S(g_j) S(g_j)^* S(g_j)
  \$=
  \varepsilon a_{kj} S(g_i)^* S(g_j),
  $$
  which can be nonzero only if $g_i=g_j$.

  \claim 3 For every $\a$ and $\b$ in $P$, if
  $S(\a)^*S(\b) \neq 0$ then $\a\inv \b \in P \cup P\inv$.

  \medskip Without loss of generality assume
  $|\a| \leq |\b|$
  and write
  $\b = \tilde \b \g$ with
  $|\tilde \b| = |\a|$ and $\tilde \b, \g \in P$.  Then
  $$
  0 \neq S(\a)^* S(\tilde \b \. \g) =
  S(\a)^* S(\tilde \b) S(\g),
  $$
  which implies that
  $S(\a)^* S(\tilde \b) \neq 0$ and hence, by claim (2), that
  $\a = \tilde \b$.  So $\a\inv \b = \g \in P \subseteq P \cup P\inv$.

  Given any $t$ in $\Fn$ let us denote by
  $E(t) = S(t) S(t)^*$.
  Since we already know that $S(\a)$ is a partial isometry for $\a\in P$ we
see that $E(\a)$ is a self-adjoint idempotent operator.

  \claim 4 For all $\a$ and $\b$ in $P$ the operators $E(\a)$ and $E(\b)$
commute.

  \medskip In the case that $\a\inv\b \not \in P \cup P\inv$ we have
  $$
  E(\a) E(\b) = S(\a) S(\a)^* S(\b) S(\b)^* = 0
  $$
  by claim (3),
  and similarly
  $E(\b) E(\a) = 0$.

  If, on the other hand ,
  $\a\inv\b \in P \cup P\inv$, without loss of generality write
  $\a\inv \b = \g \in P$ and note that
  $$
  E(\a) E(\b) =
  S(\a) S(\a)^* S(\a \. \g) S(\a \. \g)^* =
  S(\a) S(\a)^* S(\a) S(\g) S(\g)^* S(\a)^*
  \$=
  S(\a) S(\g) S(\g)^* S(\a)^* =
  S(\a) S(\g) S(\g)^* S(\a)^* S(\a) S(\a)^*
  \$=
  S(\b) S(\b)^* S(\a) S(\a)^* =
  E(\b) E(\a).
  $$

  \claim 5 For all $t$ and $r$ in $\Fn$ the operators $E(t)$ and $E(r)$
commute.

  \medskip Recall that $S(t) = 0$, unless $t=\a\b\inv$ for some $\a,\b \in P$,
which we suppose also satisfy $|t|=|\a|+|\b|$.  Then
  $S(t) = S(\a) S(\b)^*$
  and
  $$
  E(t) = S(\a) S(\b)^* S(\b) S(\a)^*.
  $$
  According to claim (1),
  $S(\b)^* S(\b) = \varepsilon S(g_i)^* S(g_i)$
  where $\varepsilon$ is either zero or one and $g_i \in G$. So
  $$
  E(t) = \varepsilon S(\a) S(g_i)^* S(g_i) S(\a)^* =
  \varepsilon \sum_{j=1}^n a_{ij} S(\a) S(g_j) S(g_j)^* S(\a)^*
  \$=
  \varepsilon \sum_{j=1}^n a_{ij} S(\a \. g_j) S(\a \. g_j)^*,
  $$
  which belongs to the linear span of the set of all $E(\g)$, with $\g\in P$;
a commutative set by claim (4).  This proves our last claim and allows us to
address the statement of the Theorem under scrutiny, that is, the proof that
  $$
  S(t) S(r) S(r\inv) = S(tr) S(r\inv)
  $$
  for $t,r \in \Fn$.

  To do this we use induction on
  $|t| + |r|$.  If either $|t|$ or $|r|$ is zero, there is nothing to prove.
So, write
  $t=\tilde t \. x$ and
  $r = y \. \tilde r$ ,
  where $x,y\in G \cup G\inv$.  In case
  $x \neq y\inv$ we have
  $|tr|=|t| + |r|$ and hence $S(t\.r) = S(t) S(r)$ as observed some time
ago. If, on the other hand,
  $x = y\inv$ we have
  $$
  S(t) S(r) S(r\inv) =
  S(\tilde t \. x) S(x\inv \. \tilde r) S(\tilde r\inv x) =
  S(\tilde t) S(x) S(x\inv) S(\tilde r) S(\tilde r\inv) S(x) = \ldots
  $$
  With an application of claim (5) and the use of the induction hypothesis we
conclude that the above equals
  $$
  \ldots =
  S(\tilde t) S(\tilde r) S(\tilde r\inv) S(x) S(x\inv) S(x) =
  S(\tilde t \tilde r) S(\tilde r\inv) S(x)
  \$=
  S(tr) S(\tilde r\inv \. x) =
  S(tr) S(r\inv).
  \proofend
  $$

  Let us briefly introduce an abstract concept, applicable to groups $\G$
which, like $\Fn$, are equipped with a ``length'' function
  $$
  |\.| \: \G \arw {\bf R}_+
  $$
  satisfying $|e| = 0$ and the triangular identity
  $|ts| \leq |t|+|s|$.

  \nstate Definition \DefSat A partial representation $\s$ of \/ $\G$ is said
to be \stress{semi-saturated} if, whenever $t,r \in \G$ are such that
  $|tr| = |t|+|r|$,
  one has
  $\s(tr) = \s(t) \s(r)$.

  Of course, the representation $S$ provided by \loccit{\PR} is a
semi-saturated representation.

  Whenever we are speaking of a partial representation $\s$ of a group $\G$,
we will denote by $\e^\s(t)$, or simply by $\e(t)$, the range projections
$\e(t)=\s(t)\s(t\inv)$, of each $\s(t)$.

  \nstate Proposition \Semisat A partial representation $\s$ of the group $\G$
is semi-saturated if and only if, whenever
  $|tr|=|t|+|r|$,
  one has that
  $\e(tr) \leq \e(t)$
  in the usual order of projections.

  \proof By using the \cstar-algebra relation $\| a\|^2 = \|aa^*\|$ and the
axioms for partial representations, one checks that, for any $t,r\in\G$,
  $$
  \| \s(t)\s(r) - \s(tr) \|^2 =
  \| \s(tr)\s(r\inv)\s(t\inv) - \e(tr) \|^2.
  $$
  Now,
  $$
  \s(tr)\s(r\inv)\s(t\inv) =
  \s(tr)\s(r\inv)\s(t\inv) \s(t) \s(t\inv)
  \$=
  \s(tr)\s(r\inv t\inv) \s(t) \s(t\inv) =
  \e(tr)\e(t),
  $$
  from which the proof can easily be concluded.
  \proofend

  Given a partial representation $S$ of $\Fn$, such as the one provided by
\loccit{\PR}, the theory developed in \cite{\Inverse} can be applied.
Specifically we would like to recall that there is a one-to-one correspondence
between partial representations of a group and representations of the partial
group \cstar-algebra, and so we get a representation of
  $C^*_p(\Fn)$.

  To be more specific, let us denote by $u_t$ the canonical partial isometry
in
  $C^*_p(\Fn)$
  associated to each $t\in \Fn$, as in the proof of \scite{\Inverse}{6.5}, and
put $\e_t = u_t u_t^*$.  Then the partial representation of $\Fn$ constructed
above gives rise, via \scite{\Inverse}{6.5}, to a representation
  $$
  \pi \: C^*_p(\Fn) \arw \Lin(H)
  $$
  such that
  $\pi(u_{g_i}) = S_i$.
  This representation is somewhat special, since the \CK relations imply that
  \smallskip \itemitem{(CK1')} $\pi(e_{g_i} e_{g_j}) = 0$, for $i\neq j$,
  \smallskip \itemitem{(CK2')} $\sum_{i=1}^n \pi(e_{g_i}) = 1$,
  \smallskip \itemitem{(CK3')} $\pi(e_{g_i\inv}) = \sum_{j=1}^n a_{i,j}
\pi(e_{g_j})$, for $i = 1, \ldots, n$.

  \medskip \noindent
  Also, by the fact that $S$ is semi-saturated, we have
  \smallskip \itemitem{(SS)} $\pi(\e_{tr}) \leq \pi(\e_t)$, whenever
$|tr|=|t|+|r|$.

  \medskip
  In other words, let $J$ be the ideal, within $C^*_p(\Fn)$, generated by the
elements
  \smallskip \itemitem{(CK1'')} $e_{g_i} e_{g_j}$, for $i\neq j$,
  \smallskip \itemitem{(CK2'')} $1 - \sum_{i=1}^n e_{g_i}$,
  \smallskip \itemitem{(CK3'')} $e_{g_i\inv} - \sum_{j=1}^n a_{i,j} e_{g_j}$
for $i = 1, \ldots, n$,
  \smallskip \itemitem{(SS')} $\e_{tr} \e_t - \e_{tr}$, for $|tr|=|t|+|r|$.

  \medskip \noindent
  Then we see that $\pi$ vanishes on $J$.  The exact logical relationship
between representations of the \CK relations and representations of
$C^*_p(\Fn)$ boils down to the following:

  \state Theorem The \CK algebra $\OA$ is canonically isomorphic to
$C^*_p(\Fn)/J$.

  \proof By \loccit{\PR} and \loccit{\Semisat} we see that there is a one to
one correspondence between representations of the \CK relations and partial
representations of $\Fn$, satisfying \loccit{CK1'--3'}+\loccit{SS}.  These, in
turn, correspond to representations of $C^*_p(\Fn)$ vanishing on $J$.  In
other words, both $\OA$ and $C^*_p(\Fn)/J$ have a universal property regarding
representations of the \CK relations and hence they must be isomorphic.
  \proofend

  Recall from \cite{\Inverse} that
  $C^*_p(\Fn)$
  is defined to be the partial crossed product of the abelian
\cstar-subalgebra $A$ generated by all the projections
  $e_t$,
  by a certain partial action of $\Fn$.  This is saying that
  $C^*_p(\Fn)$
  is the full cross-sectional algebra of the corresponding semi-direct product
bundle \cite{\TPA}.  As such, it is a topologically graded \cstar-algebra. The
grading space corresponding to the unit group element coincides with the
subalgebra $A$, mentioned above.  Incidentally, this is where the generators
of the ideal $J$ belong, and so \loccit{\Ideals} applies and tell us that the
quotient algebra, that is $\OA$, is a topologically graded algebra.  This
grading turns out to be equivalent to the co-action of $\Fn$ on $\OA$ studied
by Quigg and Raeburn in \cite{\QR}.

  By inspecting the grading in
  $C^*_p(\Fn)/J$
  from the point of view of \loccit{\Ideals}, it is easy to see that the
grading of $\OA$, thus obtained, is given by the subspaces
  $$
  {\cal O}_A^t = \overline{{\rm span}}\{ \e(r_1) \e(r_2) \ldots \e(r_k) \s(t)
\:
  k \in \N, \ r_1,r_2,\ldots,r_k \in \Fn \},
  $$
  where $\s$ is the partial representation of $\Fn$ corresponding to any given
faithful representation of $\OA$.  The Fell bundle corresponding to this
grading will henceforth be referred to as the \CK bundle.

  \section {Amenability of the \CK bundle} The main goal of this section is to
prove that the \CK bundle satisfies the approximation property.
  In fact we will obtain a somewhat more general result, as we shall see in a
moment.  The quest for wider generality will have the added advantage of
allowing the presentation to become pretty much independent of the previous
section, although, of course, the main motivation for what we are about to do
lies there.

With that in mind, let $\G$ be a discrete group (soon to be supposed equal to
$\Fn$) and let
  $$
  \sigma \: \G \arw \Lin(H)
  $$
  be any partial representation of $\G$ on $H$ and define $B_t^\s$ to be the
closed linear subspace of $\Lin(H)$ given by
  $$
  B_t^\s = \overline{{\rm span}}\{ \e(r_1) \e(r_2) \ldots \e(r_k) \s(t) \:
  k \in \N, \ r_1,r_2,\ldots,r_k \in \G\}.
  $$
  By using that
  $$
  \s(t)e(r) =
  \s(t) \s(r) \s(r\inv) =
  \s(tr) \s(r\inv)
  \$=
  \s(tr) \s(r\inv t\inv) \s(tr) \s(r\inv) =
  \s(tr) \s(r\inv t\inv) \s(t) =
  \e(tr) \s(t),
  $$
  and that
  $$
  \s(t) \s(r) =
  \s(t) \s(t\inv) \s(t) \s(r) =
  \s(t) \s(t\inv) \s(tr) =
  \e(t) \s(tr),
  $$
  one sees that
  $B_{t_1}^\s B_{t_2}^\s \subseteq B_{t_1 t_2}^\s$
  and that
  ${B_{t}^\s}^* = B_{t\inv}^\s$
  for any $t, t_1, t_2 \in \G$.

  Regardless of whether or not the $B_t^\s$ form an independent family of
linear subspaces of $\Lin(H)$, we may equip the disjoint union of the $B_t^\s$
with the Fell bundle structure over $\G$, induced by the \cstar-algebra
operations of $\B(H)$.

  \state Definition Given a partial representation of a discrete group $\G$,
we denote by $\B^\s$ the Fell bundle
  $$
  \B^\s = \(B_t^\s \)_{t\in \G}
  $$
  canonically associated to $\s$, as described above.

  Let us suppose, for the time being, that $\G$ possesses a length function,
as in our discussion before \loccit{\DefSat}.

  \state Definition A Fell bundle $\B$, over $\G$ is said to be
\stress{semi-saturated} if, whenever $t,r \in \G$ are such that
  $|tr| = |t|+|r|$,
  one has that $B_{tr}$ coincides with the closed linear span of the set of
products
  $\{b_t c_r \: b_t\in B_t, c_r \in B_r \}$ (compare
  \scite{\Circle}{4.8} and
  \scite{\FD}{VIII.2.8}).

  \state Proposition If $\s$ is a semi-saturated partial representation of \/
$\G$, then $\B^\s$ is a semi-saturated Fell bundle.

  \proof If $|t t'| = |t|+|t'|$ and if
  $r_1,r_2,\ldots,r_k, r'_1,r'_2,\ldots,r'_{k'} \in \G$ we have
  $$
  \e(r_1) \e(r_2) \ldots \e(r_k) \s(t)\quad
  \e(r'_1) \e(r'_2) \ldots \e(r'_{k'}) \s(t')
  \$=
  \e(r_1) \e(r_2) \ldots \e(r_k)
  \e(tr'_1) \e(tr'_2) \ldots \e(tr'_{k'}) \s(tt')
  $$
  which implies that $\B^\s$ is semi-saturated.
  \proofend

  In what follows we will assume that $\s$ is a semi-saturated partial
representation of $\Fn$, satisfying
  $$
  \sum_{i=1}^n \e(g_i) = 1.
  $$
  Compare with \loccit{CK2}.  We should point out that we will make no further
assumptions on $\s$ and hence we will be facing a Fell bundle that is
substantially more general than the \CK bundle.

  Observe that each $\e(g_i)$ is a projection and, because projections adding
up to one must always be orthogonal to each other, we have for $i\neq j$ that
  $$
  0 = \s(g_i)^* \e(g_i) \e(g_j) \s(g_j) =
  \s(g_i)^* \s(g_i) \s(g_i)^* \s(g_j) \s(g_j)^* \s(g_j) =
  \s(g_i)^* \s(g_j),
  $$
  telling us that the ranges of $\s(g_i)$ and $\s(g_j)$ must be orthogonal to
each other.  What this really says, is that \loccit{CK1} follows from
\loccit{CK0} and \loccit{CK2}.

  Now, since we are assuming that $\s$ is semi-saturated, this implies, as in
the case of the \CK partial representation, that a necessary condition for
$\s(t) \neq 0$ is that $t=\a\b\inv$ for some $\a,\b\in P$.

  \nstate Lemma \Soma If $k$ is a natural number then
  $ \sum_{\a \in P_k} \e(\a) = 1. $

  \proof This is obvious for $k=0$.  If $k\geq 1$ note that $P_k$ is the
disjoint union of $g_i P_{k-1}$, for $i=1,\ldots,n$. Thus, by induction we
have
  $$
  \sum_{\a \in P_k} \s(\a) \s(\a)^* =
  \sum_{i=1}^n \sum_{\b\in P_{k-1}} \s(g_i\.\b) \s(g_i\.\b)^*
  \$=
  \sum_{i=1}^n \s(g_i) \( \sum_{\b\in P_{k-1}} \s(\b) \s(\b)^* \) \s(g_i)^* =
  \sum_{i=1}^n \s(g_i) \s(g_i)^* = 1.
  \proofend
  $$

  A very important conclusion to be drawn by this is:

  \nstate Lemma \MainLemma For each $t$ in $\F_n$ and for each natural number
$k$ we have
  $$
  \s(t) = \sum_{\a\in P_k} \e(t\a) \s(t) \e(\a).
  $$

  \proof The right hand side equals
  $$
  \sum_{\a\in P_k} \s(t\a) \s(\a\inv t\inv) \s(t) \s(\a) \s(\a\inv)
  \$=
  \sum_{\a\in P_k} \s(t\a) \s(\a\inv) \s(\a) \s(\a\inv) =
  \sum_{\a\in P_k} \s(t) \s(\a) \s(\a\inv)
  = \s(t).
  \proofend
  $$

  \state Theorem For every semi-saturated partial representation $\s$ of
$\Fn$, such that
  $$
  \sum_{i=1}^n \e(g_i) = 1,
  $$
  the bundle $\B^\s$ satisfies the approximation property and hence is
amenable.

  \proof Let, for each $m\in \N$, the function
  $a_m\: \Fn \arw B_e^\s$ be defined by
  $$
  a_m(t) = \left\{ \matrix{
    m^{-{1 \over 2}}
  \e(t) &\hbox{if $t\in P$ and $|t|\leq m$}, \hfill \cr
    0 &\hbox{in all other cases}, \hfill \cr
  }\right.$$
  Observe that
  $$
  \sum_{t \in \Fn} a(t)^* a(t) =
  \sum_{k=0}^m \sum_{\a \in P_k} m\inv \e(\a)=
  \sum_{k=0}^m m\inv =
  1.
  $$

  We will show this sequence satisfies the hypothesis of \loccit{\AP}.  In
order to prove it, let us begin by showing that
  $$
  \lim_{m\rightarrow \infty} \sum_{r\in \Fn} a_m(tr)^*\s(t) a_m(r) = \s(t)
  \for t\in \Fn.
  $$

  Since $\s(t) = 0$, unless $t=\a\b\inv$ for some $\a,\b\in P$, as observed
shortly before \loccit{\Soma}, let us assume that this is the case. We may
also assume that
  $|t| = |\a| + |\b|$.

Observe that the sum above is, in fact, a finite sum, since $a_m(r)=0$ except
for $r\in P$ with $|r|\leq m$.  What is not so obvious is the behavior of
$a_m(tr)$, for these $r$'s.  Making this point clear will take us some work.

We claim that, for any $\g$ in $P$ one has
  \zitemno = 0
  \zitem if $|\g| < |\b|$ then $a_m(t\g) = 0$,
  \zitem if $|\b| \leq |\g| \leq m - |t|$ then
    $a_m(t\g) = m^{-{1 \over 2}} \e(t\g)$,

  \medskip \noindent
  To see this, note that, in the first case, $t\g$, which is given by
  $t\g = \a \. \b\inv \g$,
  cannot belong to $P$ and hence $a_m(t\g) =0 $ by definition.

  In the second case, write $\g=\g_1 \. \g_2$ with $\g_1,\g_2\in P$ and
$|\g_1|=|\b|$.  If it so happens that $\g_1 \neq \b$ then
  $t\g = \a \. \b\inv \g_1 \. \g_2 \not \in P P\inv$ which implies that
$\s(t\g)=0$ and hence also that $\e(t\g)=0$.  So, no matter which clause in
the definition of $a_m(t\g)$ is used, we will have
  $a_m(t\g) = m^{-{1 \over 2}} \e(t\g)$.

  However, if $\g_1 = \b$, we have that
  $t \g = \a \b\inv \g_1 \g_2 = \a \g_2 \in P,$
  and
  $|t \g| \leq |t| + |\g| \leq m.$
  So, by definition,
  $a_m(t\g) = m^{-{1 \over 2}} \e(t\g)$, proving our claim.

  This clarifies the role of $a_m(t\g)$ for most cases, but the case $|\g| > m
- |t|$ is still missing.  Regarding this we claim that
  \zitem If $ k > |\b|$ then
  $ \| \sum_{\g \in P_k} \e(t\g) \| \leq 1. $

  \medskip \noindent For $\g\in P_k$, write $\g=\g_1 \g_2$ with $|\g_1| =
|\b|$, repeating, somewhat, the proof of case (ii) above.  If $\g_1\neq \b$,
we saw that $\e(t\g)=0$. So, $\e(t\g)$ can only be nonzero if $\g=\b\eta$ with
$\eta\in P_{k-|\b|}$.  For $\g$ such as that we have $t\g = \a\eta$ and
therefore
  $$
  \sum_{\g \in P_k} \e(t\g) =
  \sum_{\eta \in P_{k-|\b|}} \e(\a\eta) \leq
  \sum_{\zeta \in P_{k-|\b|+|\a|}} \e(\zeta) = 1,
  $$
  the last equality following from \loccit{\Soma}.  This proves (iii).

  Now, back to the task of computing the limit referred to at the beginning of
the present proof, observe that, since $a_m(r)=0$ unless $r\in P$ and $|r|\leq
m$, we have,
  $$
  \sum_{r\in \Fn} a_m(tr)^*\s(t) a_m(r) =
  \sum_{k=0}^m \sum_{\g\in P_k} a_m(t\g)^*\s(t) a_m(\g).
  $$

  Let us now carefully analyze each sum
  $\sum_{\g\in P_k} a_m(t\g)^*\s(t) a_m(\g)$, in terms of the possible values
of $k$.
  Of course, if $k <|\b|$ this sum will vanish by (i).  Next, if
  $|\b| \leq k \leq m - |t|$, then by (ii) and \loccit{\MainLemma} we have
  $$
  \sum_{\g\in P_k} a_m(t\g)^*\s(t) a_m(\g) =
  m\inv \sum_{\g\in P_k} \e(t\g)^*\s(t) \e(\g) =
  m\inv \s(t).
  $$

  Finally, assume that
  $m - |t| < k \leq m$.  Since we are only interested in the limit as
$m\rightarrow \infty$, we may assume that $m$ is as big as necessary to force $k >
|\b|$, and hence use the conclusion of (iii).  In the computation below, we
will also use the inequality
  $$
  \|\sum_i x_i^* y_i \|^2 \leq
  \|\sum_i x_i^* x_i \| \|\sum_i y_i^* y_i \|,
  $$
  which holds in any \cstar-algebra, and can be proved by considering matrices
over the algebra.

  We have, for $k>|\b|$,
  $$
  \| \sum_{\g\in P_k} a_m(t\g)^*\s(t) a_m(\g) \|^2
  \$\leq
  \| \sum_{\g\in P_k} a_m(t\g)^*a_m(t\g) \|
    \| \sum_{\g\in P_k} a_m(\g)^*\s(t)^* \s(t) a_m(\g) \|
  \$\leq
  \| \s(t) \|^2
    \| \sum_{\g\in P_k} m\inv \e(t\g) \|
    \| \sum_{\g\in P_k} m\inv \e(\g) \|
  \leq
  m^{-2}.
  $$

  This said we have that
  $$
  \sum_{k=0}^m \sum_{\g\in P_k} a_m(t\g)^*\s(t) a_m(\g)
  \$=
  \sum_{k=|\b|}^{m-|t|} \sum_{\g\in P_k} a_m(t\g)^*\s(t) a_m(\g) +
    \sum_{k=m-|t|+1}^m \sum_{\g\in P_k} a_m(t\g)^*\s(t) a_m(\g).
  $$
  The first summand after the last equal sign equals
  $$
  \sum_{k=|\b|}^{m-|t|} m\inv\s(t) = {{m-|t|-|\b|+1}\over{m}} \s(t),
  $$
  and converges to $\s(t)$ as $m\rightarrow \infty$.

  The second one has norm no bigger than
  $$
  \sum_{k=m-|t|+1}^m \| \sum_{\g\in P_k} a_m(t\g)^*\s(t) a_m(\g) \|
  \leq
  \sum_{k=m-|t|+1}^m m\inv = {{|t|}\over{m}}
  $$
  and converges to zero as $m\rightarrow \infty$.

  This completes the task proposed at the beginning of the proof.
  Now, observe that each $B_t^\s$ is spanned by the elements of the form
  $$
  y = \e(r_1) \e(r_2) \ldots \e(r_k) \s(t),
  $$
  with
  $k \in \N$ and $r_1,r_2,\ldots,r_k \in \Fn$.  So, in view of the
commutativity of $B_e^\s$ \scite{\Inverse}{2.4} we have that
  $$
  a_m(tr)^* y a_m(r) =
  \e(r_1) \e(r_2) \ldots \e(r_k) a_m(tr)^* \s(t) a_m(r),
  $$
  and we can use the conclusion above to affirm that
  $$
  \lim_{m \rightarrow \infty} \sum_{r\in \G} a_m(tr)^* y a_m(r) = y
  $$
  for any $y$ as above.
  Finally, since the operators
  $$
  b_t \in B_t \mapsto \sum_{r\in \G} a_m(tr)^* b_t a_m(r)
  $$
  are uniformly bounded with $m$, and since the set of $y$'s as above span a
dense subset of $B_t^\s$, we conclude that
  $$
  \lim_{m \rightarrow \infty} \sum_{r\in \G} a_m(tr)^* x a_m(r) = x
  $$
  for all $x\in B_t^\s$, as required.
  \proofend

  \state Corollary The \CK bundle satisfies the approximation property and
hence is amenable.

  \proof Follows immediately from our last result.
  \proofend

  \bigbreak
  \centerline{\tensc References}
  \nobreak\medskip\frenchspacing

\ATtechreport{\AEE,
  author = {B. Abadie, S. Eilers and R. Exel},
  title = {Morita Equivalence and Crossed Products by {Hilbert} Bimodules},
  institution = {Universidade de S\~ao Paulo},
  year = {1994},
  note = {preprint},
  toappear = {Trans. Amer. Math. Soc.},
  NULL = {},
  atrib = {IR},
  MR = {}
  }

\ATarticle{\Claire,
  author = {C. Anantharaman-Delaroche},
  title = {Systemes dynamiques non commutatifs et moyenabi\-li\-t\'e},
  journal = {Math. Ann.},
  year = {1987},
  volume = {279},
  pages = {297--315},
  NULL = {},
  }

\ATarticle{\Blackadar,
  author = {B. Blackadar},
  title = {Shape theory for $C^*$-algebras},
  journal = {Math. Scand.},
  year = {1985},
  volume = {56},
  pages = {249--275},
  NULL = {},
  }

\ATarticle{\Bratteli,
  author = {O. Bratteli},
  title = {Inductive limits of finite dimensional $C^*$-algebras},
  journal = {Trans. Amer. Math. Soc.},
  year = {1972},
  volume = {171},
  pages = {195--234},
  NULL = {},
  }

\ATarticle{\CKbib,
  author = {J. Cuntz and W. Krieger},
  title = {A Class of $C^*$-algebras and Topological Markov Chains},
  journal = {Inventiones Math.},
  year = {1980},
  volume = {56},
  pages = {251--268},
  NULL = {},
  }

\ATarticle{\SoftIII,
  author = {G. Elliott, R. Exel and T. Loring},
  title = {The Soft Torus {III}: The Flip},
  journal = {J. Operator Theory},
  year = {1991},
  volume = {26},
  pages = {333--344},
  NULL = {},
  atrib = {IR},
  MR = {94f:46086}
  }

\ATarticle{\SoftI,
  author = {R. Exel},
  title = {The Soft Torus and Applications to Almost Commuting Matrices},
  journal = {Pacific J. Math.},
  year = {1993},
  volume = {160},
  pages = {207--217},
  NULL = {},
  atrib = {IR},
  MR = {94f:46091}
  }

\ATarticle{\SoftII,
  author = {R. Exel},
  title = {The Soft Torus: A Variational Analysis of Commutator Norms},
  journal = {J. Funct. Analysis},
  year = {1994},
  volume = {126},
  pages = {259--273},
  NULL = {},
  atrib = {IR},
  MR = {95i:46085}
  }

\ATarticle{\Circle,
  author = {R. Exel},
  title = {Circle Actions on {$C^*$}-Algebras, Partial Automorphisms and a
Generalized {Pimsner}--{Voiculescu} Exact Sequence},
  journal = {J. Funct. Analysis},
  year = {1994},
  volume = {122},
  pages = {361--401},
  NULL = {},
  atrib = {IR},
  MR = {95g:46122}
  }

\ATarticle{\AF,
  author = {R. Exel},
  title = {Approximately Finite {$C^*$}-Algebras and Partial Automorphisms},
  journal = {Math. Scand.},
  year = {1995},
  volume = {76},
  pages = {281--288},
  NULL = {},
  atrib = {IR},
  MR = {}
  }

\ATtechreport{\TPA,
  author = {R. Exel},
  title = {Twisted Partial Actions, A Classification of Regular
{$C^*$}-Algebraic Bundles},
  institution = {Universidade de S\~ao Paulo},
  year = {1994},
  note = {preprint},
  toappear = {Proc. London Math. Soc.},
  NULL = {},
  atrib = {N}
  }

\ATtechreport{\Unconditional,
  author = {R. Exel},
  title = {Unconditional Integrability for Dual Actions},
  institution = {Universidade de S\~ao Paulo},
  year = {1995},
  note = {preprint},
  toappear = {},
  NULL = {},
  atrib = {N}
  }

\ATtechreport{\Inverse,
  author = {R. Exel},
  title = {Partial actions of groups and actions of inverse semi-groups},
  institution = {Universidade de S\~ao Paulo},
  year = {1995},
  note = {preprint},
  toappear = {},
  NULL = {},
  atrib = {N}
  }

\ATbook{\FD,
  author = {J. M. G. Fell and R. S. Doran},
  title = {Representations of *-algebras, locally compact groups, and Banach
*-algebraic bundles},
  publisher = {Academic Press},
  year = {1988},
  volume = {125 and 126},
  series = {Pure and Applied Mathematics},
  NULL = {},
  }

\ATbook{\KK,
  author = {K. Jensen and K. Thomsen},
  title = {Elements of $K\!K$-Theory},
  publisher = {Birkh\"auser},
  year = {1991},
  volume = {},
  series = {},
  NULL = {},
  }

\ATarticle{\McClanahan,
  author = {K. McClanahan},
  title = {$K$-theory for partial crossed products by discrete groups},
  journal = {J. Funct. Analysis},
  year = {1995},
  volume = {130},
  pages = {77--117},
  NULL = {},
  }

\ATarticle{\Nica,
  author = {A. Nica},
  title = {$C^*$-algebras generated by isometries and Wiener-Hopf operators},
  journal = {J. Operator Theory},
  year = {1991},
  volume = {27},
  pages = {1--37},
  NULL = {},
  }

\ATbook{\Pedersen,
  author = {G. K. Pedersen},
  title = {$C^*$-Algebras and their automorphism groups},
  publisher = {Acad. Press},
  year = {1979},
  volume = {},
  series = {},
  NULL = {},
  }

\ATtechreport{\Quigg,
  author = {J. C. Quigg},
  title = {Discrete $C^*$-coactions and $C^*$-algebraic bundles},
  institution = {},
  year = {},
  note = {},
  toappear = {Austral. Math. J.},
  NULL = {},
  }

\ATtechreport{\QR,
  author = {J. Quigg and I. Raeburn},
  title = {Characterizations of Crossed Products by Partial Actions},
  institution = {University of Newcastle},
  year = {1996},
  note = {preprint},
  toappear = {},
  NULL = {},
  }

\ATarticle{\Rieffel,
  author = {M. A. Rieffel},
  title = {Deformation quantization of Heisenberg manifolds},
  journal = {Commun. Math. Phys.},
  year = {1989},
  volume = {122},
  pages = {531--562},
  NULL = {},
  }

\ATbook{\SV,
  author = {\c S. Str\v{a}til\v{a} and D. Voiculescu},
  title = {Representations of AF-algebras and of the group $U(\infty)$},
  publisher = {Springer-Verlag},
  year = {1975},
  volume = {486},
  series = {Lecture Notes in Mathematics},
  NULL = {},
  }

\ATbook{\Wassermann,
  author = {S. Wassermann},
  title = {Exact $C^*$-algebras and related topics},
  publisher = {Seoul National University},
  year = {1994},
  volume = {19},
  series = {Lecture Notes},
  NULL = {},
  }

\ATarticle{\Woronowicz,
  author = {S. L. Woronowicz},
  title = {Twisted $SU_2$ groups. An example of a non-commutative differential
calculus},
  journal = {Publ. RIMS, Kyoto Univ.},
  year = {1987},
  volume = {23},
  pages = {117--181},
  NULL = {},
  }

  \vskip 2cm
  \rightline{April 1996}

  \bye